\documentclass[preprint,12pt]{elsarticle} 

\usepackage[utf8]{inputenc}

\usepackage{enumerate}
\usepackage{listings}
\usepackage{longtable}
\usepackage{array}
\usepackage{graphicx}
\usepackage{enumitem}
\usepackage{amsfonts}
\usepackage{longtable}
\usepackage{float}
\usepackage{hyperref}
\usepackage{xcolor, colortbl}
\usepackage{mathtools}
\usepackage{upgreek}
\usepackage{amsmath,amsfonts,amsthm}

\usepackage[noend,noline]{algorithm2e}

\usepackage{caption}

\usepackage{subcaption}
\newtheorem{rem}{Remark}

\newcommand{\DPart}[2]{\frac{\partial #1}{\partial #2}}
\newcommand{\Vect}[1]{\mathbf{#1}}

\journal{Journal of Computational Science}

\begin{document}

\begin{frontmatter} 

\title{A study of efficient concurrent integration methods of B-Spline basis functions in IGA-FEM}

\author[label1]{Maciej Wo\'{z}niak}
\author[label1]{Anna Szyszka}
\author[label2]{Sergio Rojas}
\address[label1]{AGH University of Sciences and Technology \\ Institute of Computer Science, Electronics and Telecommunication \\ Department of Computer Science \\ al. A Mickiewicza 30, 30-059 Krak\'{o}w, Poland \\ email: macwozni@agh.edu.pl } 
\address[label2]{Instituto de Matem\'aticas, Pontificia Universidad Cat\'olica de Valpara\'iso. Valpara\'iso, Chile}

\begin{abstract} 
Based on trace theory, we study efficient methods for concurrent integration of B-spline basis functions in IGA-FEM.
We consider several scenarios of parallelization for two standard integration methods; the classical one and sum factorization.
We aim to efficiently utilize hybrid memory machines, such as modern clusters, by focusing on the non-obvious layer of the shared memory part of concurrency.
We estimate the performance of computations on a GPU and provide a strategy for performing such computations in practical implementations.
\end{abstract}


\end{frontmatter}


\section{Introduction}
\label{sec:motivation}

The great success of the finite element method (FEM) can be attributed to its solid theoretical rooting in the fields of variational calculus, and functional analysis \cite{strangfix,hughes2012finite}.
It is widely used for numerically solving partial differential equations (PDEs) in Computer-Aided Engineering (CAE) systems.
Most FEM computations consist of two phases; a) a transformation of the PDE to a discrete form by mapping onto the finite-dimensional approximation space, b) solving the resulting system of linear or nonlinear algebraic equations \cite{FEM}.
Commonly, the FEM implementation is made by computing local integral subroutines by elements, defining local element matrices that are subsequently integrated and assembled in the global system of matrix equations.

A current hot topic regarding numerical FEM approximations is the IsoGeometric Analysis FEM (IGA-FEM) \cite{Hughes}.
It integrates the geometrical modeling of CAD systems with engineering computations of CAE systems.

IGA-FEM computations share the same structure as the traditional FEM.
However, the main difference is that IGA-FEM employs B-splines basis functions for spanning the approximation space \cite{SubroutinePackageForCalculating}.

In several scenarios, mainly when dealing with time-dependency or Non-linearity, it is well known that FEM can give rise to the resolution of a high-cost computational problem.
For instance, one of the standard techniques for numerically solving time-dependent PDEs is to perform a finite difference method (FDM) in time, coupled with a FEM discretization in space.
Last implies assembling multiple FEM matrices at every time step in several scenarios.
Indeed,  for nonlinear PDEs, it may be required to integrate and assemble FEM matrices at each iteration step of the nonlinear solver.
In particular, if an implicit method in time is employed \cite{isotumor3d, Puzyrev, CHICCS2016}.
Furthermore, the cost of the assembling grows with the space dimension \cite{cost}.
Therefore, the cost associated with the integration in assembling FEM matrices is critical in terms of computation time.

Traditionally, the integration procedure is performed in parallel, element-by-element, making a level of concurrent operations in which data is independent.
However, in \cite{parallel_integration} is proposed a methodology based on adding two levels of parallelism within each element that distinguishes the independent operations.
The goal was to reduce the computational time in the integration procedure using the modern parallel architectures of a GPU \cite{CUDA}.

The paper aims to compare the practical concurrent implementation performance of the classical integration method and sum factorization with different parallelization schemes.
To do that, we apply the methodology presented in \cite{parallel_integration} to sum factorization.
For this, we start by briefly describing the principal concepts involved in the case study.

\subsection{Architecture}
\label{sec:architecture}

\begin{figure}[ht]
    \centering
    \includegraphics[width=0.95\textwidth]{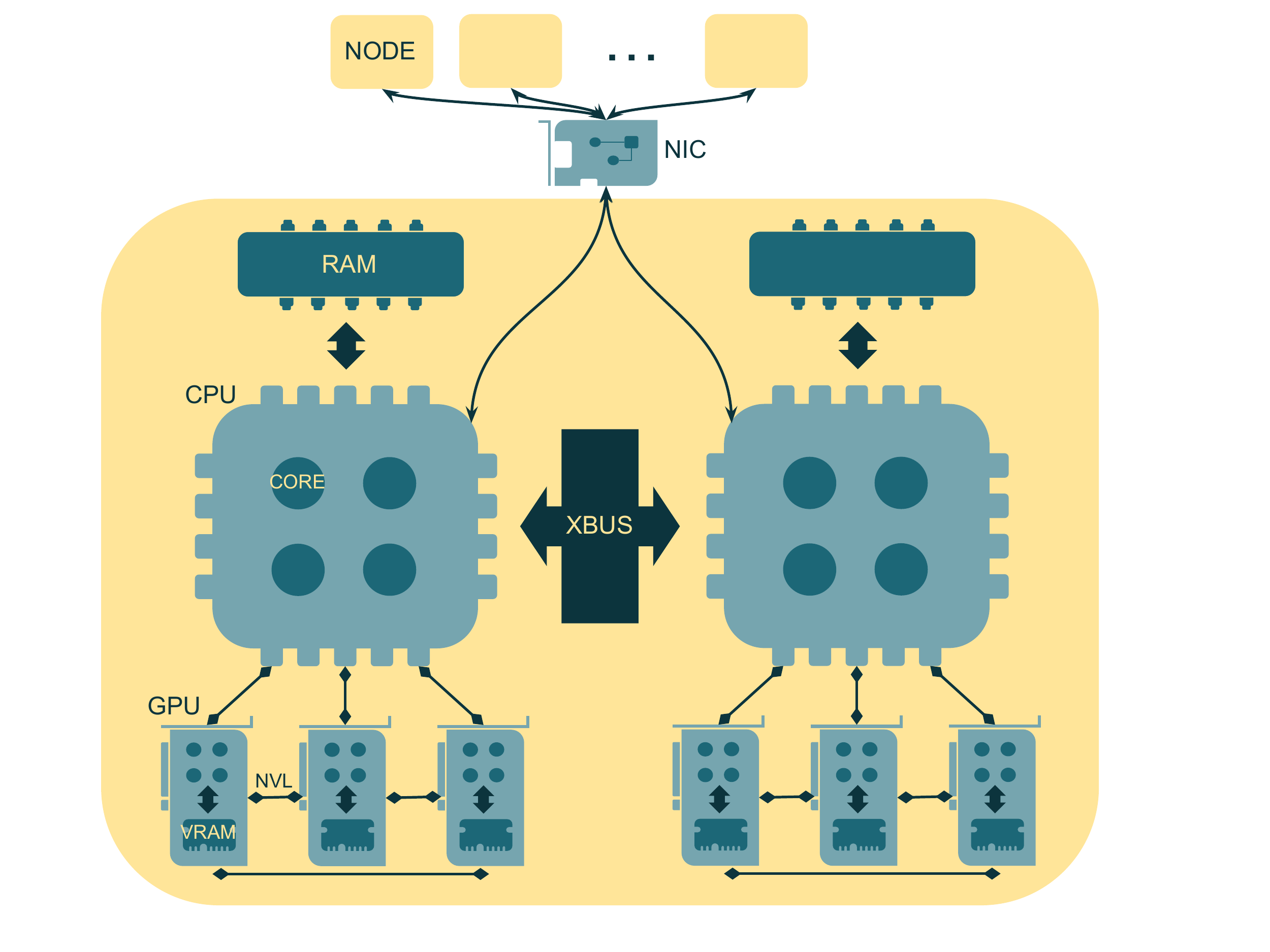}
    \caption{Architecture of a single node of a modern supercomputer.
    Figure based on Summit \cite{Summit} supercomputer.}
    \label{fig:node}
\end{figure}

State-of-the-art supercomputers are designed as multi-level hierarchical hybrid systems \cite{cyfronet,stampede,Summit}.
A representative architecture is shown in Figure~\ref{fig:node}.
They consist of classical nodes (servers) communicating over a network (specialized solutions, such as Infiniband) through a Message Passing Interface (MPI).
Inside every server (node), multiple multi-core CPUs partially share RAM.
Furthermore, these systems have a massively parallel co-processor such as GPUs.
GPUs have dedicated memory, with a hierarchical memory organization \cite{CUDAmemory}, which is not shared with the CPU.
Concurrent algorithms dedicated to such systems are crucial for efficient hardware utilization and reduced carbon trace (green computing).

\subsection{Sum factorization}
\label{sec:sumfact_intro}

Sum factorization (see, e.g.~\cite{REV1}) was first introduced in \cite{sumfactbook}.
It was initially employed for the standard higher-order finite element method \cite{REV0}.
However, currently it is preferred the technique of choice for efficient formation of local element matrices in hp-finite elements \cite{VOS20105161, 10.1137/11082539X, karniadakis2005spectral, eibner2005fast} and IGA with higher-order B-splines \cite{REV1,ANTOLIN2015817, BRESSAN2019437}.
In essence, sum factorization is a reordering of the computations in such a way as to exploit the underlying tensor product of the test and trial spaces involved.
When employed, the cost of integration is reduced from ${\cal O}(m^3n^3q^3)$ to ${\cal O}(m^3n^3q+m^2n^2q+mnq^3)$, where $m$ denotes the number of test functions over an element in each direction, $n$ is the number of trial functions over the element, and $q$ is the number of quadrature points over the element in each spatial direction.
Last in practice implies that, for a given polynomial degree $p$, the total reduction is from ${\cal O}(p^9)$ to ${\cal O}(p^7)$ when considering Gaussian quadrature, and down up to ${\cal O}(p^6)$ for weighted quadrature.

\subsection{Trace theory}
\label{sec:tracetheory}

There exist multiple methods used in the formal verification of concurrent computations.
One of the most popular is the \textit{Trace Theory} \cite{TraceTheory}.
Other methods include \textit{Petri Net} \cite{PetriNet}, \textit{Process Calculi} \cite{ProcessCalculi} and \textit{Actor Model} \cite{ActorModel}.

Trace Theory delivers the Foata Normal Form (FNF) \cite{BookOfTraces} and Diekert dependency graphs, which help characterize the processing in a single element and simplify the parallel implementation on massively parallel machines, such as GPUs.
Finally, it makes a base for near-optimal scheduling. It simplifies concurrent implementation on a GPU while providing a theoretical framework for verifying the correctness of such a parallel algorithm.

\subsection{Structure of the article}
The rest of the article is organized as follows.
First, in Section \ref{sec:model} we describe the model problem, together with its discretization in time and space, used for the benchmarks.
Next, in Section~\ref{sec:integration} we discuss the integration algorithms and apply trace theory to create a concurrent algorithm performing sum factorization.
In Section~\ref{sec:numres} we consider several numerical experiments to show and discuss the 
performance of the integration methodologies.
Finally, we conclude the paper in Section \ref{sec:conclussions}.

\section{Model problem and IGA-discrete variational formulation}
\label{sec:model}

\subsection{Model problem}
\label{sec:formulation}

With the spirit of presenting the proposed methodology in a simple setting (i.e., the extension \cite{parallel_integration}
to the concurrent sum factorization algorithm), 
we will consider the following heat-transfer model problem:   
\begin{equation}
  \left\{
    \begin{aligned}
      \DPart{u}{t} &=
      \Delta u 
      &\qquad&\text{in }\Omega\times[0, T]\\
      \nabla u \cdot \hat{\Vect{n}} &= 0
      &\qquad&\text{on }\partial\,\Omega\times[0,T] \\
      u &= u_0 \text{, at } t=0 &\qquad&\text{in }\Omega
    \end{aligned}
  \right.
  \label{eq:heat}
\end{equation}
where $\Omega = (0 , 1)^3  \subset \mathbb{R}^3$ denotes the spatial domain, $\hat{\Vect{n}}$ denotes the normal vector to the domain boundary $\partial \Omega$, $T>0$~is a length of the time interval, and $u_0$ is a given initial state.

\subsection{Discretization in time}
To obtain a fully-discrete formulation of problem~\eqref{eq:heat}, we start by considering its corresponding continuous weak formulation in space, given as follows: 

Find~$u \in \mathcal{C}^1\left(\left(0, T\right), H^1\left(\Omega\right)\right)$ such that $u=u_0$ at $t=0$ and, for each~$t \in \left(0, T\right)$, it holds:
\begin{equation}
  \label{eq:weak_space}
  \int_\Omega \DPart{u}{t} v \, dx =
  - \int_\Omega \nabla u \cdot \nabla v \, dx, \, \qquad \forall \, v \in H^1\left(\Omega\right).
\end{equation}
For simplicity, we consider a discrete-in-time version of problem~\eqref{eq:weak_space} by employing the forward Euler method.
This is, denoting by $u_n$ the approximation of $u$ at time $t = n\Delta_t$, with $n=0,\dots,N$ where $\Delta_t = T/N$ denotes a fixed time step for a given integer $N>0$, we obtain $u_{n+1} \in H^1(\Omega)$ as the solution of the following variational problem:
\begin{equation}
  \label{eq:Euler}
  \int_\Omega u_{n+1} v \, dx = 
   \int_\Omega u_n v \, dx - \Delta_t \int_\Omega \nabla u_n \cdot \nabla v \, dx,  \qquad \forall \, v \in H^1\left(\Omega\right).
\end{equation}

\subsection{B-splines basis functions}
\label{sec:splines}

A B-spline is a convenient function representing polynomial splines (see e.g.~\cite{boor,shumaker}).
B-splines are characterized by the polynomial degree inside the respective elements and their regularity at the interfaces between them of the finite element mesh.
For simplicity, in this work, we will consider 3D-tensor B-splines basis functions of the same polynomial degree and regularity at the interior faces of the tensor mesh.
However, the methodology can be easily extended to more general B-spline basis functions.

Consider a partitioning of $\overline{\omega} = [ 0 , 1 ]$ into $K$ uniform elements $[\widehat{x}_{k-1},\widehat{x}_k]$, with
$$ \widehat{x}_0 = 0 < \widehat{x}_1 < \dots < \widehat{x}_{k-1} < \widehat{x}_{k} < \dots < \widehat{x}_{K} = 1.$$
For a given $p>0$, the B-spline basis functions being piece-wise polynomials of degree $p$, with $C^{p-1}$ regularity at the interior knots $\{\widehat{x}_k\}_{k=1}^{K-1}$, are defined trough the following knot vector:
\begin{equation}\label{eq:knot_vector}
    \Xi = \{ \gamma_i \}^{K+p-1}_{i=0} := 
    \{ \underbrace{0,\dots,0}_{p+1} , \dots , \widehat{x}_{k-1}, \widehat{x}_k, \dots, \underbrace{1,\dots,1}_{p+1} \}.
\end{equation}
More precisely, the $i$-th B-spline basis function, with $1\leq i\leq K + p - 1$, is constructed using the Cox--de--Boor recursive formulae \cite{SubroutinePackageForCalculating}:
\begin{equation}
    B_{i;0}(\xi):=\left\{\begin{array}{ll} 1, & \mathrm {if} \quad \gamma_{i}\leq \xi <\gamma_{i+1} , \\ 0, & \mathrm {otherwise} ,
    \end{array}\right.
\label{eq:cox1}
\end{equation}
\begin{equation}
    B_{i;q}(\xi):={\frac {\xi-\gamma_{i}}{\gamma_{i+q}-\gamma_{i}}}B_{i;q-1}(\xi)+{\frac {\gamma_{i+q+1}-\xi}{\gamma_{i+q+1}-\gamma_{i+1}}}B_{i+1;q-1}(\xi), \text{ for } 1\leq q \leq p,
\label{eq:cox2}
\end{equation}
where $B_{i,q}(\xi)$ denotes the value of the $i$-th B-spline function of degree $q$ at the point $\xi$.
In formula (\ref{eq:cox2}), the limit case $0/0$ is defined as $0$.
We notice that the Cox--de--Boor recursive formulae \eqref{eq:cox1} and \eqref{eq:cox2} define a total of $K+p$ 1D B-splines basis functions.\\
We define 3D basis functions by tensor product of 1D B-splines basis functions that, for simplicity, we construct considering the same number of element partitions for all three spatial directions.\\
For $x = (x_1,x_2,x_3) \in \Omega$, we will denote by 
\begin{equation}\label{eq:basis_construction}
B_{\delta;p}(x) = B_{k;p}(x_1) B_{l;p}(x_2) B_{m;p}(x_3), \text{ with } \delta = \{k,l,m\} \in \mathcal{K},
\end{equation}
the evaluation in $x$ of a generic 3D B-spline basis function, where 
\begin{equation}
\mathcal{K} = \{0,1,\dots,K+p-1\}^3.
\end{equation}
Finally, we will denote by 
\begin{equation}\label{eq:B_spline_space}
\mathcal{B}_{\mathcal{K};p}:=\text{span}\left\{B_{\delta;p}, \text{ with } \delta \in \mathcal{K}\right\} \subset H^1(\Omega)
\end{equation}
the space generated by the 3D--tensor B--spline basis functions of degree $p$ and global regularity $p-1$.

\subsection{Fully-discrete variational formulation}
For a given polynomial degree $p$, the fully--discrete formulation of problem~\eqref{eq:heat} is obtained from \eqref{eq:Euler} by considering the $H^1$--conforming space $\mathcal{B}_{\mathcal{K};p}$ as the approximation space for the discrete solution $U_{n+1} \approx u_{n+1}$.
This is, given $U_{n} \in \mathcal{B}_{\mathcal{K};p}$, we obtain $U_{n+1} \in \mathcal{B}_{\mathcal{K};p}$ as the solution of the following discrete variational formulation problem:
\begin{equation}\label{eq:discrete_variational_formulation}
\textrm{Find } U_{n+1} \in \mathcal{B}_{\mathcal{K};p}, \textrm{ such that } a\left(U_{n+1},B_{\delta;p}\right)=l_n\left(B_{\delta;p}\right), \forall \,\delta \in \mathcal{K},
\end{equation}
where
\begin{equation}
a\left(U_{n+1},B_{\delta;p}\right)=\int_{\Omega} U_{n+1} B_{\delta;p}\, dx,
\end{equation}
\begin{equation}
l_n\left(B_{\delta;p}\right)=\int_\Omega U_{n} B_{\delta;p}\, dx - \Delta t \int_{\Omega} \nabla U_{n} \cdot \nabla B_{\delta;p} \, dx,
\end{equation}
and $U_{0}$ corresponds to the classical $L^2$-projection of the initial state $u_0$ in the B-spline space $\mathcal{B}_{\mathcal{K};p}$.\\
As a consequence of the finite number of basis functions for the discrete space $\mathcal{B}_{\mathcal{K};p}$, we can assume that the wanted discrete solution is written as:
\begin{equation}
U_{n+1}(x) = \sum_{\beta \in \{1, \dots, K+p\}^3} \mu_\beta B_{\beta;p}(x).
\end{equation}
Therefore, after considering an appropriate ordering for the basis functions that here we will consider implicit for the sake of simplicity, problem~\eqref{eq:discrete_variational_formulation} can be equivalently written in matrix form as:
\begin{equation}\label{eq:matrix_form}
\text{Find } \mu \in \mathbb{R}^{(K+p)^3}, \text{ such that } A \mu = L,
\end{equation}
with the right-hand side $L_\delta = l_n \left(B_{\delta;p}\right)$, and the Gram matrix
\begin{equation}\label{eq:mass_term}
A_{\delta,\beta} = a\left(B_{\beta;p}, B_{\delta;p}\right).
\end{equation}

\section{Integration algorithms}
\label{sec:integration}

\subsection{Element-by-element integration}
For the sake of simplicity, here we will focus on the integration and assembling of the Gram matrix $A$.

The standard integration strategy consists of assembling the linear system \eqref{eq:matrix_form} element-by-element.  
To exemplify the procedure, we assume that the domain $\Omega$ is decomposed into a set of $K^3$ cubic elements.
\begin{equation}\label{eq:Edelta}
   E_{\gamma} = (\widehat{x}_i,\widehat{x}_{i+1}) \times (\widehat{x}_j,\widehat{x}_{j+1}) \times (\widehat{x}_k,\widehat{x}_{k+1}), 
\end{equation}
where $\widehat{x}_j = j/K$ (cf.~Section~\ref{sec:splines}), and $\gamma = (i,j,k) \in  \{1, 2, \ldots, K\}^3$.

Denoting by $\delta = (h,i,j)$, and by $\beta = (k,l,m)$, the matrix element $A_{\delta,\beta}$ (see \eqref{eq:mass_term}) is computed as the sum 
\begin{equation}
A_{\delta,\beta} = \sum_{\gamma \in \{1, 2, \ldots, K\}^3} A_{\delta,\beta}^{\gamma},
\end{equation}
where $A_{\delta, \beta}^{\gamma}$ is given in terms of the 1D B-spline basis functions as (see \eqref{eq:basis_construction}):
\begin{equation}
A_{\delta, \beta}^{\gamma} = \int_{E_{\gamma}}
B_{h;\, p}(x_1) \, B_{i;\, p}(x_2) \, B_{j;\, p}(x_3) \, B_{k;\, p}(x_1) \, B_{l;\, p}(x_2) \, B_{m;\, p}(x_3) \, dx.
\label{eq:88}
\end{equation}
Let us consider a proper exact quadrature with the particular set of weights and nodes
$\{\omega^n = (\omega^{n_1},\omega^{n_2}, \omega^{n_3})$,
$x^n  = (x^{n_1}, x^{n_2}, x^{n_3}) \in E_{\gamma}\}$, with
$n_1=1,\dots,P_1$, 
$n_2=1,\dots,P_2$,
$n_3=1,\dots,P_3$,
$n=1,\ldots,P$, 
and $P = P_1 P_2 P_3$ depending on the quadrature rule and polynomial order $p$. Then, the matrix element \eqref{eq:88} is computed as:
\begin{equation}
\displaystyle A_{\delta, \beta}^{\gamma} = \sum_{n_1=1}^{P_1} \sum_{n_2=1}^{P_2} \sum_{n_3=1}^{P_3}  \omega^{n_1} \omega^{n_2} \omega^{n_3} \, \Pi(x^n) \, J(x^n) \, dx, 
\label{eq:10}
\end{equation}
where $\Pi(x^n) = B_{h;\, p}(x^{n_1}) \, B_{i;\, p}(x^{n_2}) \, B_{j;\, p}(x^{n_3}) \,
B_{k;\, p}(x^{n_1}) \, B_{l;\, p}(x^{n_2}) \, B_{m;\, p}(x^{n_3}) $ and $J(x^n)$ corresponds to the Jacobian of the particular element evaluated at $x^n$.

\begin{rem}[Element-by-element pre-computations]
We notice that~\eqref{eq:10} can be efficiently calculated by first pre-computing, over each element, only the integral of the B-splines products with non-empty support. Therefore, for a given $\alpha=(i,j,k)$, it will be helpful to introduce the set of multi-indices
\begin{equation}
\mathcal{K}^\Delta_{\alpha} = \{ (z,r,s): z \in \{i, \dots, i+p\}, r \in \{j, \dots, j+p\}, s \in \{k, \dots, k+p\}\}
\end{equation}
corresponding to the indexes of the $(p+1)^3$ basis functions with non-empty support in the $\alpha$-element.
\end{rem}

\subsection{Algorithm descriptions and computational cost}
\label{sec:algorithm}

In this section, we describe the two algorithms to be compared in subsequent sections, the classical integration algorithm, and the sum factorization algorithm.

On one side, in the {\bf classical integration} algorithm, local contributions to the left-hand-side Gram matrix are represented as a sum over quadrature points, as shown in equation (\ref{eq:10}) and described in Algorithm~\ref{algorithm1}.
In this case, the associated computational cost is known that scales, concerning the polynomial degree $p$, as ${\cal O}(p^9)$ \cite{HIEMSTRA2019234}.
\begin{algorithm}
\For{element $E \in \Omega $}{
\For{test function $B_{i,x}$}{
\For{trial function $B_{j,x}$}{
\For{test function $B_{i,y}$}{
\For{trial function $B_{j,y}$}{
\For{test function $B_{i,z}$}{
\For{trial function $B_{j,z}$}
{
    \For{quadrature point $(\xi, w)$ in $E$}{
      $A_{\delta,\beta} \gets A_{\delta,\beta} + B_{i,x}(\xi_x)B_{j,x}(\xi_x) \, B_{i,y}(\xi_y) B_{j,y}(\xi_y) \, B_{i,z}(\xi_z) B_{j,z}(\xi_z) \, J(\xi) \, w$\;
    }
  }}}
}}}
}
\caption{Classical integration algorithm}
\label{algorithm1}
\end{algorithm}

On another side, {\bf Sum factorization} algorithm consists of reorganizing the integration terms of equation \eqref{eq:10} to reduce the computational cost, in terms of the polynomial degree $p$, associated with the sum procedure.
In practice, equation \eqref{eq:10} is written as:
\begin{equation}
  A_{\beta, \delta} = \sum_{n_3=1}^{P_3} \omega^n_3 \, B_{j;\, p}(x^n_3) \, B_{m;\, p}(x^n_3) \,
   \textcolor{red}{C(i_2,i_3,j_2,j_3,k_1)},
   \label{eq:finalsum}
\end{equation}
where buffer~$C$ is given by
\begin{equation}
  \textcolor{red}{C(i_2,i_3,j_2,j_3,k_1)} = 
  \sum_{n_2=1}^{P_2} \omega^n_2 \, B_{i;\, p}(x^n_2) \, B_{l;\, p}(x^n_2)
  \textcolor{blue}{\underbrace{\sum_{n_1=1}^{P_1} \omega^n_1 \, B_{h;\, p}(x^n_1) \, B_{k;\, p}(x^n_1) \, J(x^n)}_{D(i_3,j_3,k_1,k_2)}}.
  \label{eq:bufferC}
\end{equation}
The algorithm is described in Algorithm~\ref{algorithm2}.
Here we can observe three distinct groups of loops.
As a consequence, this implies that the computational cost associated with sum factorization is ${\cal O}(p^7)$ \cite{HIEMSTRA2019234}.
\begin{algorithm}

\For{test function $B_{i,z}$ - ($i_3$)}{
\For{trial function $B_{j,z}$ - ($j_3$)}{
\For{quadrature point $(\xi_x, w_x)$ in $E$ - ($k_1$)}{
\For{quadrature point $(\xi_y, w_y)$ in $E$ - ($k_2$)}{
  \For{quadrature point $(\xi_z, w_z)$ in $E$ - ($k_3$)}{
    $D(i_3,j_3,k_1,k_2) \gets D(i_3,j_3,k_1,k_2) + B_{i,z}(\xi_z) B_{j,z}(\xi_z) \, w_z \, J(\xi)$\;
  }
}}
}}
\For{test function $B_{i,y}$ - ($i_2$)}{
\For{trial function $B_{j,y}$ - ($j_2$)}{
\For{test function $B_{i,z}$ - ($i_3$)}{
\For{trial function $B_{j,z}$ - ($j_3$)}{
  \For{quadrature point $(\xi_x, w_x)$ in $E$ - ($k_1$)}{
  \For{quadrature point $(\xi_y, w_y)$ in $E$ - ($k_2$)}{
      $C(i_2,i_3,j_2,j_3,k_1) \gets C(i_2,i_3,j_2,j_3,k_1) + B_{i,y}(\xi_y) B_{j,y}(\xi_y) \, D(i_3,j_3,k_1,k_2) \, w_y$\;
  }}
}}
}}
\For{test function $B_{i,x}$ - ($i_1$)}{
\For{trial function $B_{j,x}$ - ($j_1$)}{
\For{test function $B_{i,y}$ - ($i_2$)}{
\For{trial function $B_{j,y}$ - ($j_2$)}{
\For{test function $B_{i,z}$ - ($i_3$)}{
\For{trial function $B_{j,z}$ - ($j_3$)}{
  \For{quadrature point $(\xi_x, w_x)$ in $E$ - ($k_1$)}{
      $A(i_1,j_1) \gets A(i_1,j_1) + B_{i,x}(\xi_x)B_{j,x}(\xi_x) \, C(i_2,i_3,j_2,j_3,k_1) \, w_x$\; \smallskip
    }
  }}
}}
}}
\caption{Sum factorization algorithm}
\label{algorithm2}
\end{algorithm}

\subsection{Concurrency model for sum factorization}
\label{sec:concurencymodel}

Multiple methods are used to verify concurrent computations by creating a concurrency model formally.
In~\cite{parallel_integration}, a concurrency model based on the \textit{Trace Theory}, introduced by Diekert and Mazurkiewicz \cite{TraceTheory}, is discussed.
It contains four levels of concurrency:
\begin{enumerate}
    \item Concurrent computations on parts of the mesh.
    \item Concurrent computations on single elements.
    \item Concurrent computations of single entries in an element matrix.
    \item Concurrent computations of Cox--de--Boor formulae and 3D B-spline functions evaluation.
\end{enumerate}
Using the same methodology, we will discuss the last two levels of concurrency for the sum factorization algorithm.

\noindent
The alphabet of tasks for the integration of B-Spline basis functions over a given element consists of the following nine tasks:
\begin{enumerate}
    \item 
    $t_{\alpha;d}^{0;r;n}$
    - computational task evaluating a 1D basis function with subscript $r$ and order $0$, over the element $E_{\alpha}$ at the coordinate of quadrature point $x_d^n$.
    Task $t_{\alpha;d}^{0;r;n}$ refers to formula \eqref{eq:cox1}.
    Namely, it computes the function $B_{r;0}(x_{d}^{n})$ over the element $E_{\alpha}$.
    
    \item $t_{\alpha;d}^{p;r;n}$, ($p > 0$)
    - computational task evaluating a 1D basis functions with subscript $r$ and order $p$, over the element $E_{\alpha}$ at the coordinate of quadrature point $x_d^n$.
    Task $t_{\alpha;d}^{p;r;x}$ refers to formula \eqref{eq:cox2}.
    It contains a series of sums, subtractions, multiplications, and divisions, using output from tasks $t_{\alpha;d}^{p-1;r;n}$ and $t_{\alpha;d}^{p-1;r+1;n}$.
    Namely, it computes the function $B_{r;p}(x_{d}^{n})$ over the element $E_{\alpha}$.
    
    \item $s_{\alpha}^{p;n}$ 
    - computational task evaluating the Jacobian value $J(x^n)$ over the element $E_\alpha$.
    Namely, it computes $J(x^n)$ according to formula (\ref{eq:10}).
    
    \item $t_{\alpha;1}^{p;\beta,\gamma;n}$
    - computational task evaluating the value of the product of two 1D basis functions $B^1_{\beta;p}$, $B^1_{\gamma;p}$, and the Jacobian value $J(x^n)$, over the element $E_{\alpha}$, at the quadrature point $x^n$.
    Task $t_{\alpha;1}^{p;\beta,\gamma;n}$ consists of a multiplication of output from tasks  $t_{\alpha;1}^{p;{a};n}$, $t_{\alpha;1}^{p;{b};n}$, and $s_{\alpha}^{p;n}$.
    Namely, it computes buffer $K_{\beta,\gamma;p}(x^n) = B_{a;p}(x_1^n)B_{b;p}(x_1^n)J(x^n)$ according to formula (\ref{eq:bufferC}).
    
    \item $s_{\alpha;1}^{p;\beta,\gamma;n}$
    - computational task evaluating the sum of  $K_{\beta,\gamma;p}(x^n)$ along $x_1$.
    Task $s_{\alpha;1}^{p;\beta,\gamma;x}$ consists of a sum of outputs from tasks $t_{\alpha;1}^{p;\beta,\gamma;n}$.
    Namely, it computes the buffer $C_{\beta,\gamma;p}(x^n) = \sum\limits_{n=1}^{P_1} K_{\beta,\gamma;p}(x^n)$ according to formula (\ref{eq:bufferC}).
    
    \item $t_{\alpha;2}^{p;\beta,\gamma;n}$
    - computational task evaluating the value of product of two 1D basis function $B_{\beta;p}$, $B_{\gamma;p}$, and sums it with the previous buffer value $C_{\beta,\gamma;p}$, over the element $E_{\alpha}$, at the quadrature point $x^n_2$.
    Task $t_{\alpha;2}^{p;\beta,\gamma;x}$ consists of a multiplication of output from tasks  $t_{\alpha;2}^{p;{a};n}$, $t_{\alpha;2}^{p;{b};n}$, and $t_{\alpha;1}^{p;\beta,\gamma;n}$.
    Namely, it computes the buffer $F_{\beta,\gamma;p}(x^n) = B_{a;p}(x_2^n)B_{b;p}(x_2^n) + C_{\beta,\gamma;p}(x^n)$ according to formula (\ref{eq:bufferC}).
    
    \item $s_{\alpha;2}^{p;\beta,\gamma;n}$
    - computational task evaluating the sum of  $E_{\beta,\gamma;p}(x^n)$ along $x_1$.
    Task $s_{\alpha;2}^{p;\beta,\gamma;x}$ consists of a sum of outputs from tasks $t_{\alpha;2}^{p;\beta,\gamma;n}$.
    Namely, it computes the buffer $D_{\beta,\gamma;p}(x^n) = \sum\limits_{n=1}^{P_2} F_{\beta,\gamma;p}(x^n)$ according to formula (\ref{eq:bufferC}).
    
    \item $t_{\alpha;3}^{p;\beta,\gamma;n}$
    - computational task evaluating the value of the product of two 1D basis functions $B_{\beta;p}$, $B_{\gamma;p}$, and sums it with the previous buffer value $D_{\beta,\gamma;p}$, over element $E_{\alpha}$, at the quadrature point $x^n_3$.
    Task $t_{\alpha;3}^{p;\beta,\gamma;x}$ consists of a multiplication of outputs from tasks  $t_{\alpha;3}^{p;{a};n}$, $t_{\alpha;3}^{p;{b};n}$, and $t_{\alpha;2}^{p;\beta,\gamma;n}$.
    Namely, it computes $H_{\beta,\gamma}^\alpha = B_{a;p}(x_3^n)B_{b;p}(x_3^n) + D_{\beta,\gamma;p}(x^n)$ according to formula (\ref{eq:finalsum}).
    
    \item $s_{\alpha;3}^{p;\beta,\gamma;n}$
    - computational task evaluating the sum of  $E_{\beta,\gamma;p}(x^n)$ along $x_1$.
    Task $s_{\alpha;3}^{p;\beta,\gamma;x}$ consists of a sum of outputs from tasks  $t_{\alpha;1}^{p;\beta,\gamma;n}$.
    Namely, it computes the buffer $A_{\beta,\gamma} = \sum\limits_{n=1}^{P_3} H_{\beta,\gamma;p}(x^n)$ according to formula (\ref{eq:finalsum}).
    
\end{enumerate}

Summarizing, each task has two, three, or four upper subscripts in two or three groups divided by a semicolon.
The first group $p$ determines the B-spline order.
The second group (optional) of multi indexes $\beta$, $\gamma$ or index $r$ determines B-spline functions indexes.
The third (optional) group $n$ determines quadrature point $x^n$, at which the functions are evaluated.
Additionally, tasks have one or two bottom subscripts.
The first one, with the index $\alpha$, determines the element over which we perform computations.
The second (optional) determines the direction in the $x$, $y$, or $z$ axis.
It is important to recall that a particular task cannot be performed until the completion of the tasks for which its output is required.

\subsubsection{Set of dependencies}
\label{sec:setdepen}

\begin{figure}[htp!]
    \centering
    \includegraphics[width=1.0\textwidth]{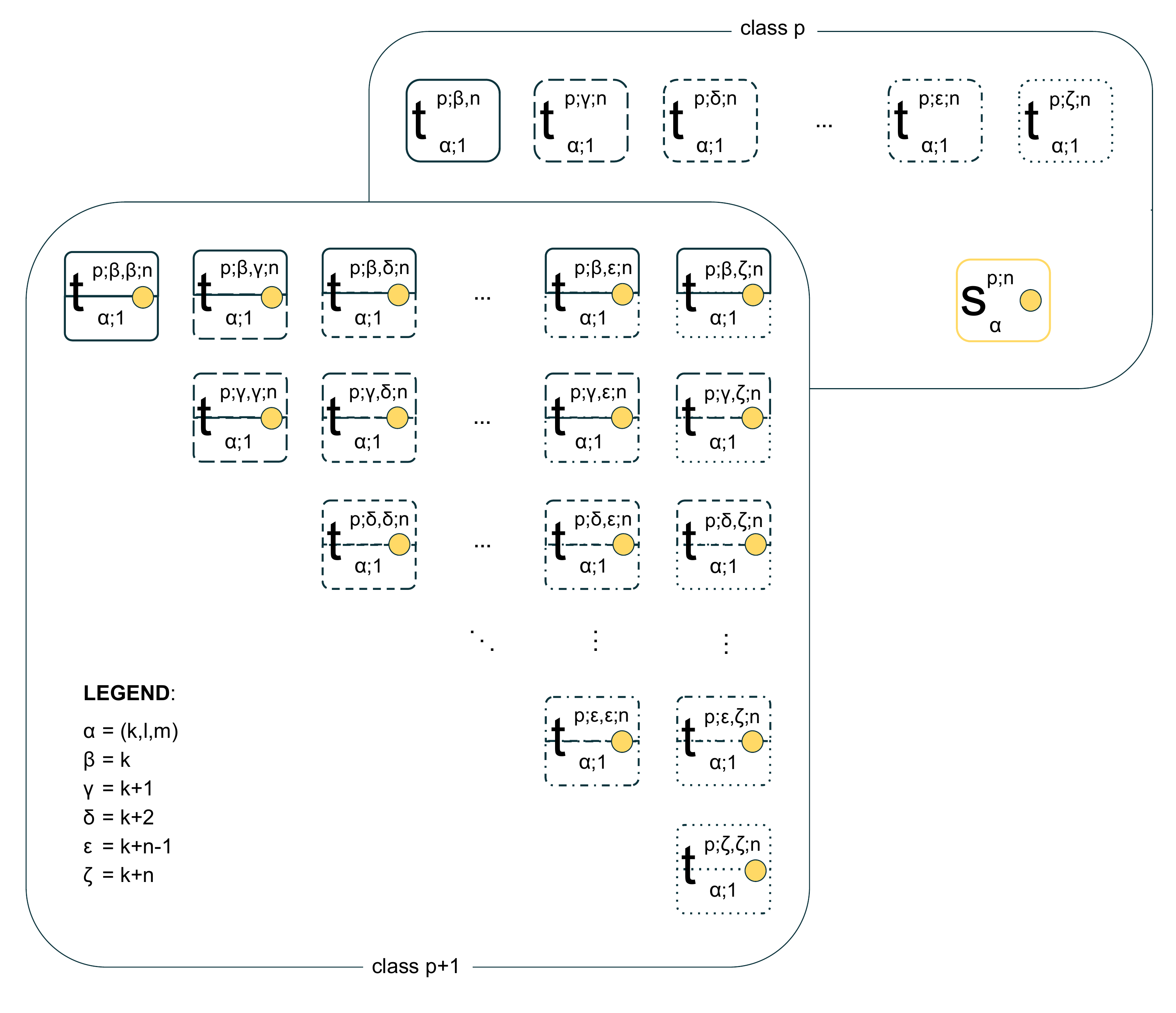}
    \caption{Relationships between classes $p$ and $p+1$ for $p$ - order functions.
    Each task in class $p+1$ corresponds to the dot product of two 1D B-spline functions, so it depends on the two tasks in class $p$.
    To maintain the transparency of the chart, the relationships between the second and third classes are marked with a border type.}
    \label{fig:n_part2}
\end{figure}

In this section, we define the alphabet of tasks $\Sigma$ and the set of dependencies between them denoted by $D$.
For this, we start by setting the variables:
\begin{align}
    & n \in \{1,2,\dots,P\}, \nonumber \\
    & r \in \{0,1,\dots,p \}, \nonumber \\
    & d \in \{1,2,3\} \nonumber. \\
    & f \in \{ k, k+1, \dots, k+p \}, \nonumber \\
    & g \in \{ l, l+1, \dots, l+p \}, \\
    & h \in \{ m, m+1, \dots, m+p \}, \nonumber \\
    & \alpha = (k,l,m) \in \{1,\dots,K\}^3, \nonumber \\
    & \beta = (a,b,c) \in \mathcal{K}^\Delta_\alpha, \nonumber \\
    & \gamma \in \mathcal{K}^\Delta_\alpha. \nonumber
    \label{eq:ranges}
\end{align}
We also set $I_\alpha$ as the function computing the index in the local element ($E_\alpha$) matrix based on the multi index $(a,b,c)$.
This is, 
\begin{equation}
I_\alpha : \mathcal{K}_\alpha^\Delta \rightarrow \{0,1,2,\dots,(p+1)^3 \}.
\label{eq:index_function}
\end{equation}
We define the alphabet of tasks as:
\begin{equation}
\Sigma = \left\{ t_{\alpha;1}^{r;f;n} , t_{\alpha;2}^{r;g;n}, t_{\alpha;3}^{r;h;n},
t_{\alpha}^{p;\beta;n}, s^{p;n}_{\alpha} \right\}
\cup \left\{
t_{\alpha;d}^{p;\beta,\gamma;n} , s_{\alpha;d}^{p;\beta,\gamma}; \, I_\alpha(\beta) \geq I_\alpha(\gamma) \right\},
\label{eq:alphabet}
\end{equation}
and the set of dependencies between tasks from the alphabet $\Sigma$ as:
\begin{align}
D = & \, J^+ \cup (J^+)^{-1}  \cup I_\Sigma, \label{eq:formulaD}
\end{align}
where
\begin{align}
J = & \, J_1 \cup J_2 \cup J_3 \cup J_4, 
\label{eq:formulaJ}
\end{align}
with
\begin{align*}
J_1 = & \, \Big\{ ( t_{\alpha;1}^{r-1;f;n} , t_{\alpha;1}^{r;f;n} ) , ( t_{\alpha;1}^{r-1;f+1;n} , t_{\alpha;1}^{r;f;n} ) , ( t_{\alpha;2}^{r-1;g;n} , t_{\alpha;2}^{r;g;n} ) , \nonumber\\
& \vspace{2cm}( t_{\alpha;2}^{r-1;g+1;n} , t_{\alpha;2}^{r;g;n} ), ( t_{\alpha;3}^{r-1;h;n} , t_{\alpha;3}^{r;h;n} ) , ( t_{\alpha;3}^{r-1;h+1;n} , t_{\alpha;3}^{r;h;n} )\Big\}, \\
J_2 = & \left\{ ( t_{\alpha;1}^{p;\alpha;n} , t_{\alpha;1}^{p;\beta,\gamma;n}) , ( t_{\alpha;2}^{p;\alpha;n} , t_{\alpha;2}^{p;\beta,\gamma;n}) , ( t_{\alpha;3}^{p;\alpha;n} , t_{\alpha;3}^{p;\beta,\gamma;n}) \right\}, \\
J_3 = & \left\{ (t_{\alpha;1}^{p;\beta,\gamma;n}, s_{\alpha;1}^{p;\beta,\gamma;n}), (t_{\alpha;2}^{p;\beta,\gamma;n}, s_{\alpha;2}^{p;\beta,\gamma;n}) , (t_{\alpha;3}^{p;\beta,\gamma;n}, s_{\alpha;3}^{p;\beta,\gamma;n}) \right\}, \\
J_4 = & \left\{ (s_{\alpha}^{p;n}, t_{\alpha;1}^{p;\beta,\gamma;n}) , (s_{\alpha;1}^{p;\beta,\gamma;n}, t_{\alpha;2}^{p;\beta,\gamma;n}) , (s_{\alpha;2}^{p;\beta,\gamma;n}, t_{\alpha;3}^{p;\beta,\gamma;n}) \right\}. 
\end{align*}
Primitives described above define the monoid of traces for the problems under consideration.
$J$ defined in equation~\eqref{eq:formulaJ} will stand for edges in Diekert dependency graph\cite{TraceTheory}, which will be drawn later in frame of this model in Figures \ref{fig:n_part1}-\ref{fig:n_part7}.

After building the primitives of the trace monoid. This is, the alphabet of tasks (\ref{eq:alphabet}) and the dependency relation (\ref{eq:formulaD}),
we define the pseudo-code allowing to compute the value of integral (\ref{eq:88}), presented in Tables \ref{tab:algorithm1}-\ref{tab:algorithm3}, that we have split into three parts to facilitate its reading.
The dependencies in this algorithm's record determine only the sequence of operations in one string representing the desired trace.
The alphabet of tasks $\Sigma$ (\ref{eq:alphabet}), the dependencies relation $D$ (\ref{eq:formulaD}), and the trace defined by pseudocode (Tables \ref{tab:algorithm1}-\ref{tab:algorithm3}) allow us to compute the Diekert dependency graph, which is convenient for the correct and effective scheduling of tasks in a heterogeneous computer environment.

\lstset
{ 
    basicstyle=\fontsize{9}{11}\selectfont\ttfamily,
    numbers=left,
    stepnumber=1,
    showstringspaces=false,
    tabsize=1,
    breaklines=true,
    breakatwhitespace=false,
}

\begin{table}[htp!]
\centering
\begin{tabular}{c}
\noindent\begin{minipage}{1.05\linewidth}
\begin{lstlisting}[escapeinside={(*}{*)}, frame = single, framexleftmargin=20pt]
(*\textbf{BEGIN}*)
//loop over elements
(*\textbf{FOREACH}*) (*$\{\alpha::=(k,l,m)\} \in \mathcal{K}^\Delta$*)
    //compute local element matrix
    element_matrix = zeros(*$\left((p+1)^3,(p+1)^3\right)$*)
    local_matrix = zeros(*$\left((p+1)^3,(p+1)^3,P\right)$*)
    local_C_matrix = zeros((*$P_y,P_z,p+1,p+1,P_x$*))
    element_C_matrix = zeros((*$P_y,P_z,p+1,p+1$*))
    local_D_matrix = zeros((*$P_z,p+1,p+1,p+1,p+1,P_x$*))
    element_D_matrix = zeros((*$P_z,p+1,p+1,p+1,p+1$*))
    
    //loop over quadrature points
    (*\textbf{FOR}*) (*$n_x$*)=1,(*$P_x$*)
        1D_matrix = zeros((*$p+1$*))
        //compute 1D functions
        (*\textbf{FOR}*) (*$r$*)=0,(*$p$*)
            (*$t_{\alpha;1}^{p;k;n}$*): 1D_matrix((*$r$*)) = compute recursive (* $\left( B_{k+r;p}(x^n_1) \right)$ *)
        (*\textbf{ENDFOR}*)
        (*\textbf{FOR}*) (*$n_y$*)=1,(*$P_y$*)
        (*\textbf{FOR}*) (*$n_z$*)=1,(*$P_z$*)
            (*$s_{\alpha}^{p;n}$*): (*$c=J(x^n)$*)
            //compute product of two functions
            (*\textbf{FOREACH}*) (*$ \beta = (a,:,:) \in \mathcal{K}^\Delta_\alpha $*)
                (*$i$*) = index_in_local_matrix((*$a,:,:$*))
                (*\textbf{FOREACH}*) (*$ \gamma = (d,:,:) \in \mathcal{K}^\Delta_\alpha$*)
                    (*$j$*) = index_in_local_matrix((*$d,:,:$*))
                    (*$t_{\alpha;1}^{p;\beta,\gamma;n}$*): local_C_matrix((*$n_y,n_z,i,j,n_x$*)) =
                    = 1D_matrix((*$a$*)) * 1D_matrix((*$d$*)) * c
                (*\textbf{ENDFOR}*)
            (*\textbf{ENDFOR}*)
        (*\textbf{ENDFOR}*)
        (*\textbf{ENDFOR}*)
    (*\textbf{ENDFOR}*)
    //sum local components from each quadrature point
    (*\textbf{FOREACH}*) (*$ \beta = (a,:,:) \in \mathcal{K}^\Delta_\alpha $*)
        (*$i$*) = index_in_local_matrix((*$a,:,:$*))
        (*\textbf{FOREACH}*) (*$ \gamma = (d,:,:) \in \mathcal{K}^\Delta_\alpha$*)
            (*$j$*) = index_in_local_matrix((*$d,:,:$*))
            (*$s_{\alpha;1}^{p;\beta,\gamma}$*): element_C_matrix((*$n_y,n_z,i,j$*)) =
            = reduction(local_C_matrix((*$n_y,n_z,i,j,:$*)),+)
        (*\textbf{ENDFOR}*)
    (*\textbf{ENDFOR}*)
    
\end{lstlisting}
\end{minipage}
\end{tabular}
\caption{The algorithm generating sample string of tasks representing sum factorization for the Gram matrix.
Part 1.}
\label{tab:algorithm1}
\end{table}
\begin{table}[htp!]
\centering
\begin{tabular}{c}
\noindent\begin{minipage}{1.05\linewidth}
\begin{lstlisting}[escapeinside={(*}{*)}, frame = single, framexleftmargin=20pt]
    //loop over quadrature points
    (*\textbf{FOR}*) (*$n_y$*)=1,(*$P_y$*)
        1D_matrix = zeros((*$p+1$*))
        //compute 1D functions
        (*\textbf{FOR}*) (*$r$*)=0,(*$p$*)
            (*$t_{\alpha;2}^{p;k;n}$*): 1D_matrix((*$r$*)) = comformulaJ3pute recursive (* $\left( B_{k+r;p}(x^n_2) \right)$ *)
        (*\textbf{ENDFOR}*)
        (*\textbf{FOR}*) (*$n_z$*)=1,(*$P_z$*)
            //compute product of two functions
            (*\textbf{FOREACH}*) (*$ \beta = (a,b,:) \in \mathcal{K}^\Delta_\alpha $*)
                (*$[i_1,i_2]$*) = index_in_local_matrix((*$a,b,:$*))
                (*\textbf{FOREACH}*) (*$ \gamma = (d,e,:) \in \mathcal{K}^\Delta_\alpha$*)
                    (*$[j_1,j_2]$*) = index_in_local_matrix((*$d,e,:$*))
                    (*$t_{\alpha;2}^{p;\beta,\gamma;n}$*): local_D_matrix((*$n_z,i_1,i_2,j_1,j_2,n_y$*)) =
                    = 1D_matrix((*$b$*)) * 1D_matrix((*$e$*)) *
                    * element_C_matrix((*$n_y,n_z,i_1,j_1$*))
                (*\textbf{ENDFOR}*)
            (*\textbf{ENDFOR}*)
        (*\textbf{ENDFOR}*)
    (*\textbf{ENDFOR}*)
    //sum local components from each quadrature point
    (*\textbf{FOREACH}*) (*$ \beta = (a,b,:) \in \mathcal{K}^\Delta_\alpha $*)
        (*$[i_1,i_2]$*) = index_in_local_matrix((*$a,b,:$*))
        (*\textbf{FOREACH}*) (*$ \gamma = (d,e,:) \in \mathcal{K}^\Delta_\alpha$*)
            (*$[j_1,j_2]$*) = index_in_local_matrix((*$d,e,:$*))
            (*$s_{\alpha;2}^{p;\beta,\gamma}$*): element_D_matrix((*$n_z,i_1,i_2,j_1,j_2$*)) =
            = reduction(local_D_matrix((*$n_z,i_1,i_2,j_1,j_2,:$*)),+)
        (*\textbf{ENDFOR}*)
    (*\textbf{ENDFOR}*)
    
\end{lstlisting}
\end{minipage}
\end{tabular}
\caption{The algorithm generating sample string of tasks representing sum factorization for the Gram matrix.
Part 2.}
\label{tab:algorithm2}
\end{table}
\begin{table}[htp!]
\centering
\begin{tabular}{c}
\noindent\begin{minipage}{1.05\linewidth}
\begin{lstlisting}[escapeinside={(*}{*)}, frame = single, framexleftmargin=20pt]
    //loop over quadrature points
    (*\textbf{FOR}*) (*$n_z$*)=1,(*$P_z$*)
        1D_matrix = zeros((*$p+1$*))
        //compute 1D functions
        (*\textbf{FOR}*) (*$r$*)=0,(*$p$*)
            (*$t_{\alpha;3}^{p;k;n}$*): 1D_matrix((*$r$*)) = compute recursive (* $\left( B_{k+r;p}(x^n_3) \right)$ *)
        (*\textbf{ENDFOR}*)
        //compute product of two functions
        (*\textbf{FOREACH}*) (*$ \beta = (a,b,c) \in \mathcal{K}^\Delta_\alpha $*)
            (*$[i;i_1,i_2,i_3]$*) = index_in_local_matrix((*$a,b,c$*))
            (*\textbf{FOREACH}*) (*$ \gamma = (d,e,f) \in \mathcal{K}^\Delta_\alpha$*)
                (*$[j;j_1,j_2,j_3]$*) = index_in_local_matrix((*$d,e,f$*))
                (*$t_{\alpha;3}^{p;\beta,\gamma;n}$*): local_matrix((*$i,j,n_z$*)) =
                = 1D_matrix((*$c$*)) * 1D_matrix((*$f$*)) *
                * element_D_matrix((*$n_z,i_1,i_2,j_1,j_2$*))
            (*\textbf{ENDFOR}*)
        (*\textbf{ENDFOR}*)
    (*\textbf{ENDFOR}*)
    //sum local components from each quadrature point
    (*\textbf{FOREACH}*) (*$ \beta = (a,b,c) \in \mathcal{K}^\Delta_\alpha $*)
        (*$[i;i_1,i_2,i_3]$*) = index_in_local_matrix((*$a,b,c$*))
        (*\textbf{FOREACH}*) (*$ \gamma = (d,e,f) \in \mathcal{K}^\Delta_\alpha$*)
            (*$[j;j_1,j_2,j_3]$*) = index_in_local_matrix((*$d,e,f$*))
            (*$s_{\alpha;3}^{p;\beta,\gamma}$*): element_matrix((*$i,j$*)) =
            = reduction(local_matrix((*$i,j,n_z$*)),+)
        (*\textbf{ENDFOR}*)
    (*\textbf{ENDFOR}*)
    
    
    //insert local matrices into global ones
    insert_local_element_2_global(element_matrix,(*$\alpha$*))
(*\textbf{ENDFOR}*)
(*\textbf{END}*)
\end{lstlisting}
\end{minipage}
\end{tabular}
\caption{The algorithm generating sample string of tasks representing sum factorization for the Gram matrix.
Part 3.}
\label{tab:algorithm3}
\end{table}

\subsection{Application of trace theory to sum factorization}
\label{sec:porder}

This section describes the methodology for creating the Diekert Dependency Graph (DG) and the Foata Normal Form (FNF), applied to the sum factorization integration method of $p$-order B-spline basis functions.
DG presents all computational tasks performed in computation and dependencies between them.
Within DG and FNF, we can distinguish Foata classes, which help with practically implementing concurrent computations.

For a given polynomial degree $p$, there are $(p+1)^3$ basis functions with non-empty support over each cubic element $E_\alpha$, with $\alpha \in \mathcal{K}^\Delta$.
Therefore, for every $E_\alpha$, we require to construct a Gram element matrix of size $(p+1)^3\times(p+1)^3$, according to equation \eqref{eq:88}.
However, due to the symmetry of the Gram matrix, it is not necessary to compute the full element matrix.
Indeed, we only require to compute $(p + 1 + (p+1)^3\times(p+1)^3)/2$ matrix entries.

To exemplify the cost associated with the computation of a single entry in the Gram matrix,  let us assume that a quadrature of $P=P_xP_yP_z$ points per element is employed.
Let us also denote by $x^1,x^2,\dots,x^P$ the corresponding quadrature points, where $x^n = (x_1^n, x_2^n, x_3^n)$, for $n=1,\dots,P$.

\begin{figure}[ht!]
    \centering
    \includegraphics[width=0.95\textwidth]{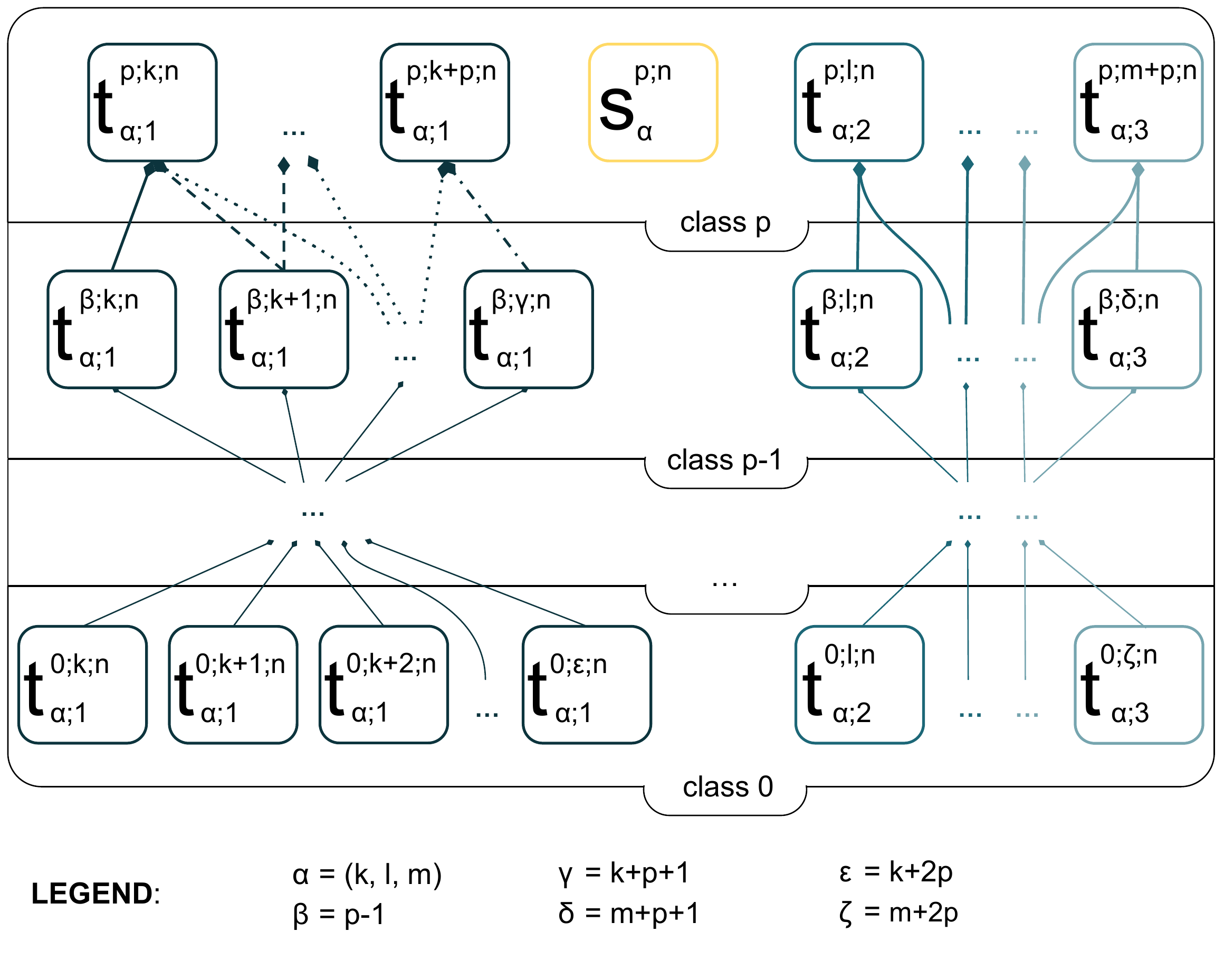}
    \caption{Relationships between classes $0$ to $p$ for $p$ - order functions.
    Tasks belonging to one class correspond to going through one iteration of the Cox-de Boor recursion formulae (\ref{eq:cox1}, \ref{eq:cox2}) for each of the three dimensions of the model.
    The dimensions are differentiated by color.}
    \label{fig:n_part1}
\end{figure}

In the procedure for each quadrature point, we start by computing $(p+1)$ 1D functions in each direction ($3p+3$ functions in total) employing the Cox--de--Boor formulae (Classes $0,1,\dots,p$ in \mbox{Figure \ref{fig:n_part1}}).
This completes all tasks of type $t_{\alpha,d}^{0;r;n}$, $t_{\alpha,d}^{1;r;n}$ up to $t_{\alpha,d}^{p;r;n}$ (see Table \ref{table66}).
Within the class $p$, we include one extra task computing $s_{\alpha}^{p;n}$.

\begin{figure}[ht!]
    \centering
    \includegraphics[width=1.0\textwidth]{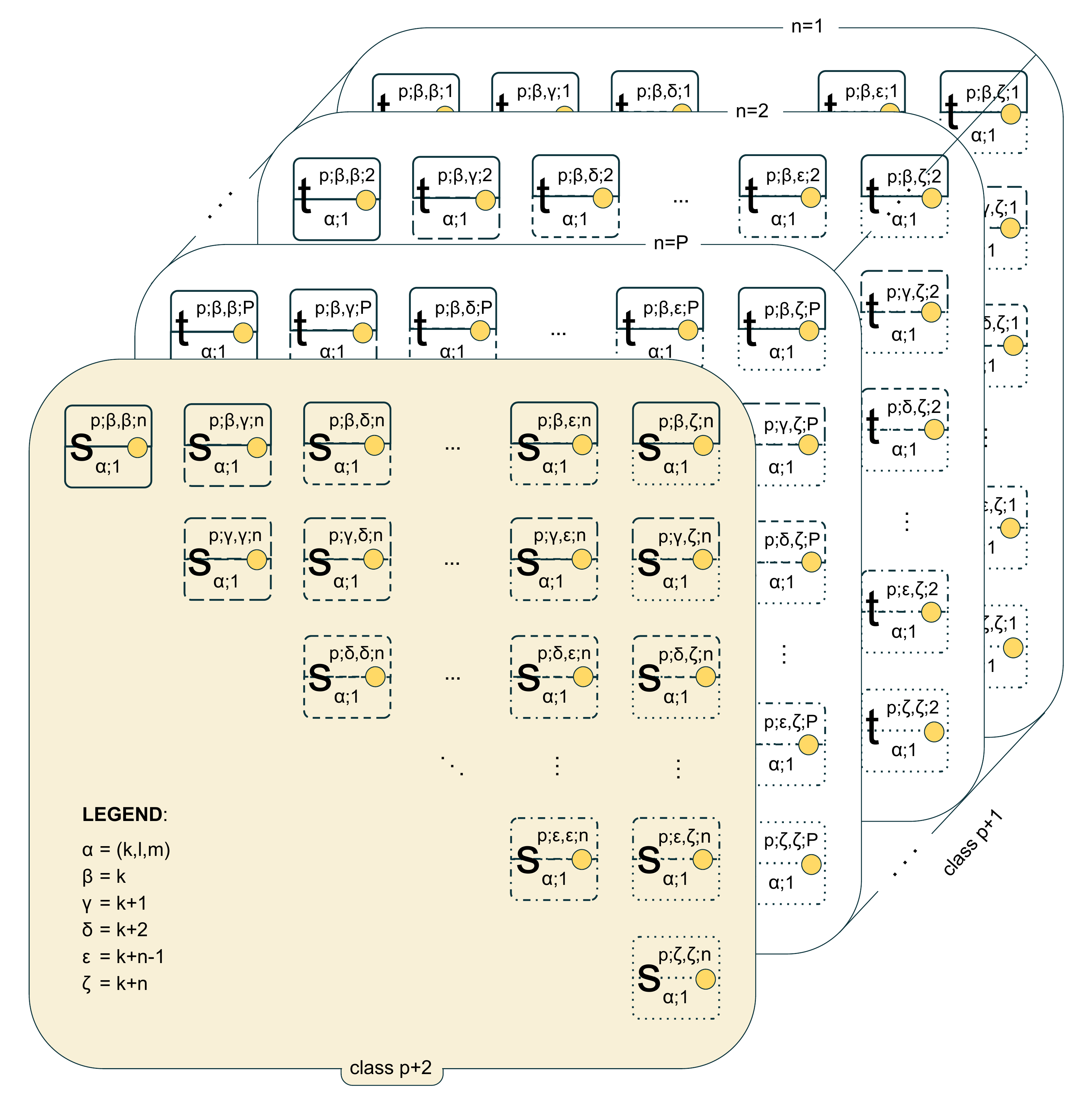}
    \caption{Relationships between classes $p+1$ and $p+2$ for $p$ - order functions.
    Each task from class $p+3$ corresponds to the approximation of the function value using Gaussian quadrature, therefore it depends on $P$ tasks from class $p+1$.
    The task $s$ depends on all tasks $t$ with regards to its distribution in subsequent sheets $1, 2, \dots, P$.}
    \label{fig:n_part3}
\end{figure}

The class $p+1$, in \mbox{Figure \ref{fig:n_part3}}), completes all tasks of the type $t_{\alpha;1}^{p;\beta,\gamma;n}$ (see Table \ref{table66}).

The concurrently computed components can be summed to evaluate scalar products of the 1D basis functions over the element $E_{\alpha}$, $(k,l,m)=\alpha \in \mathcal{K}^\Delta$, which completes all tasks of the type $s_{\alpha;1}^{p;\beta,\gamma}$ (see Table \ref{table66}).

\begin{figure}[htp!]
    \centering
    \includegraphics[width=1.0\textwidth]{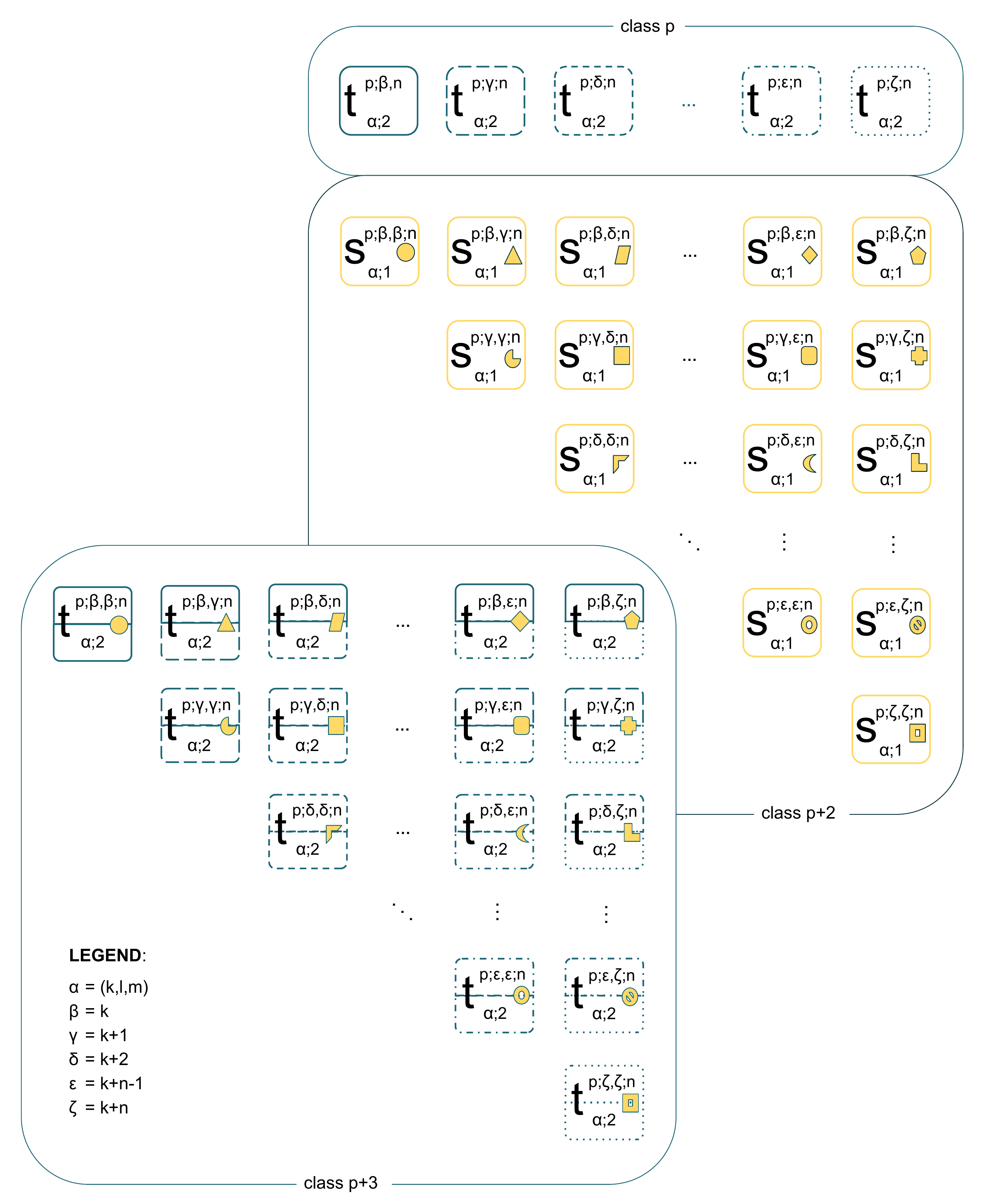}
    \caption{Relationships between classes $p$, $p+2$, and $p+3$ for $p$ - order functions.
    Each task in class $p+1$ corresponds to the dot product of two 1D B-spline functions and a sum, so it depends on the two tasks in class $p$ and some tasks from class $p+2$.
    To maintain the transparency of the chart, the relationships between the second and third classes are marked with a border type.}
    \label{fig:n_part4}
\end{figure}

\begin{figure}[ht!]
    \centering
    \includegraphics[width=1.0\textwidth]{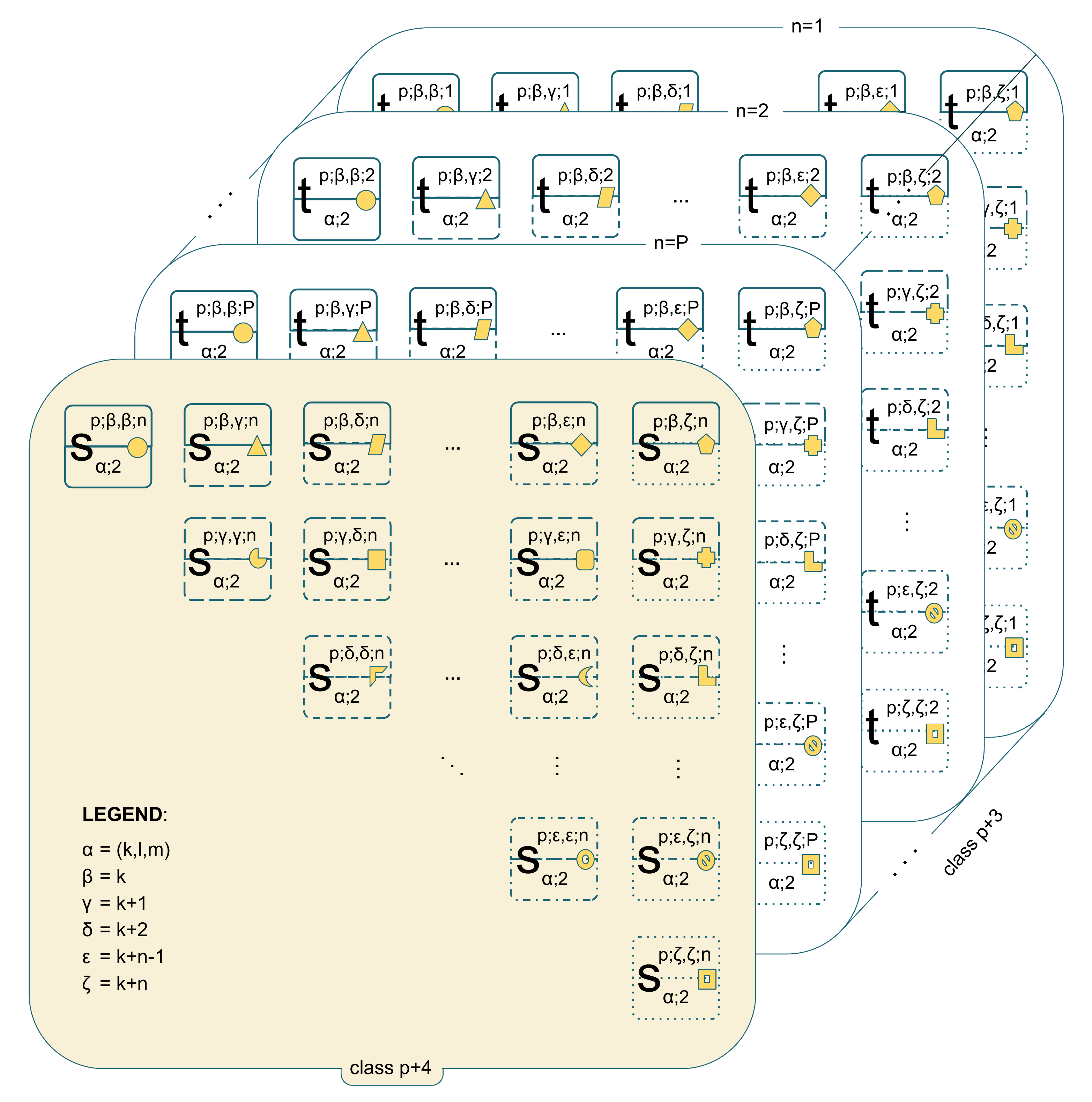}
    \caption{Relationships between classes $p+3$ and $p+4$ for $p$ - order functions.
    Each task from class $p+3$ corresponds to the approximation of the function value using Gaussian quadrature, therefore it depends on $P$ tasks from class $p+3$.
    The task $s$ depends on all tasks $t$ with regards to its distribution in subsequent sheets $1, 2, \dots, P$.}
    \label{fig:n_part5}
\end{figure}

\begin{figure}[htp!]
    \centering
    \includegraphics[width=1.0\textwidth]{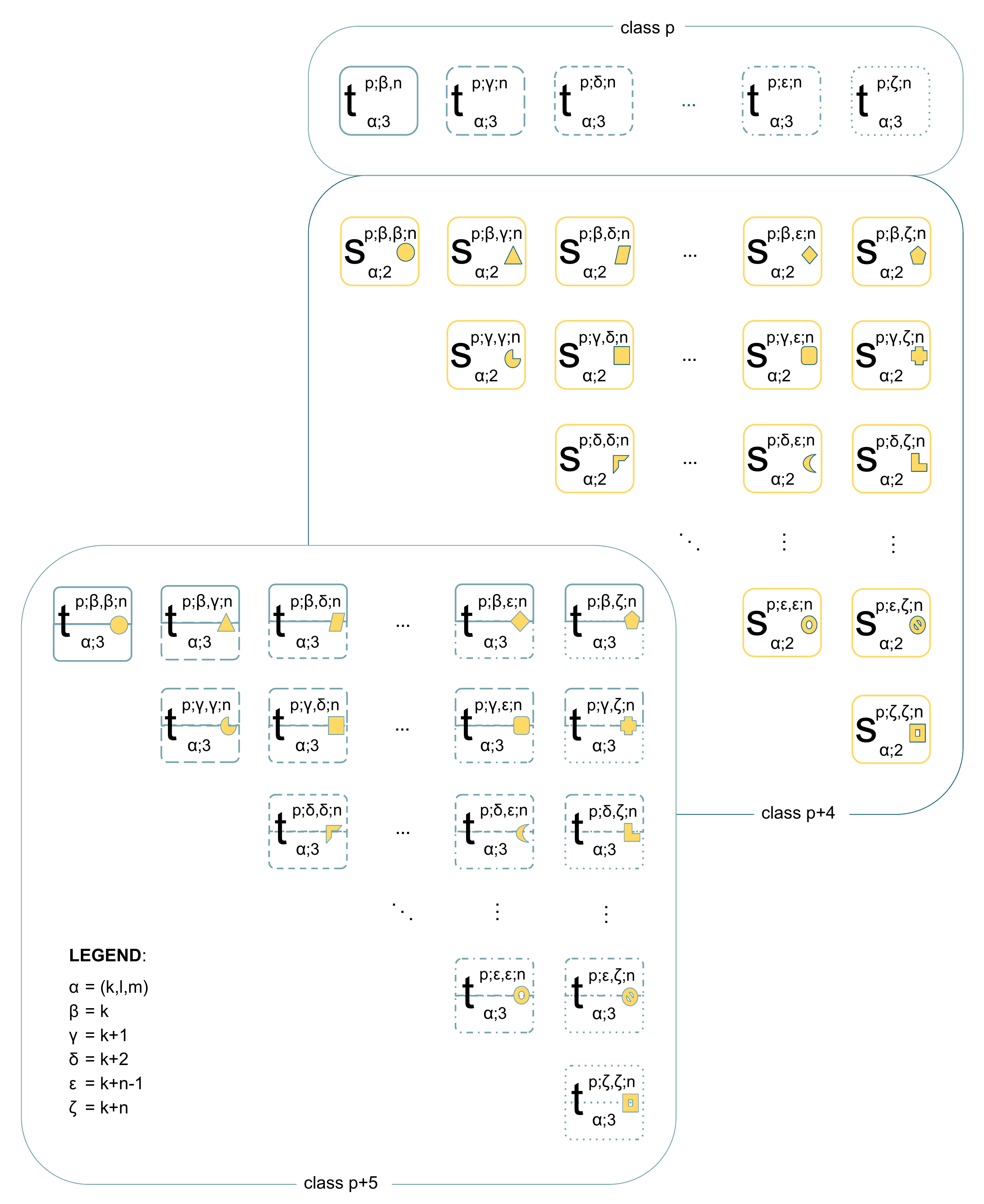}
    \caption{Relationships between classes $p$, $p+4$, and $p+5$ for $p$ - order functions.
    Each task in class $p+5$ corresponds to the dot product of two 1D B-spline functions and a sum, so it depends on the two tasks in class $p$ and some tasks from class $p+4$.
    To maintain the transparency of the chart, the relationships between the second and third classes are marked with a border type.}
    \label{fig:n_part6}
\end{figure}

\begin{figure}[ht!]
    \centering
    \includegraphics[width=1.0\textwidth]{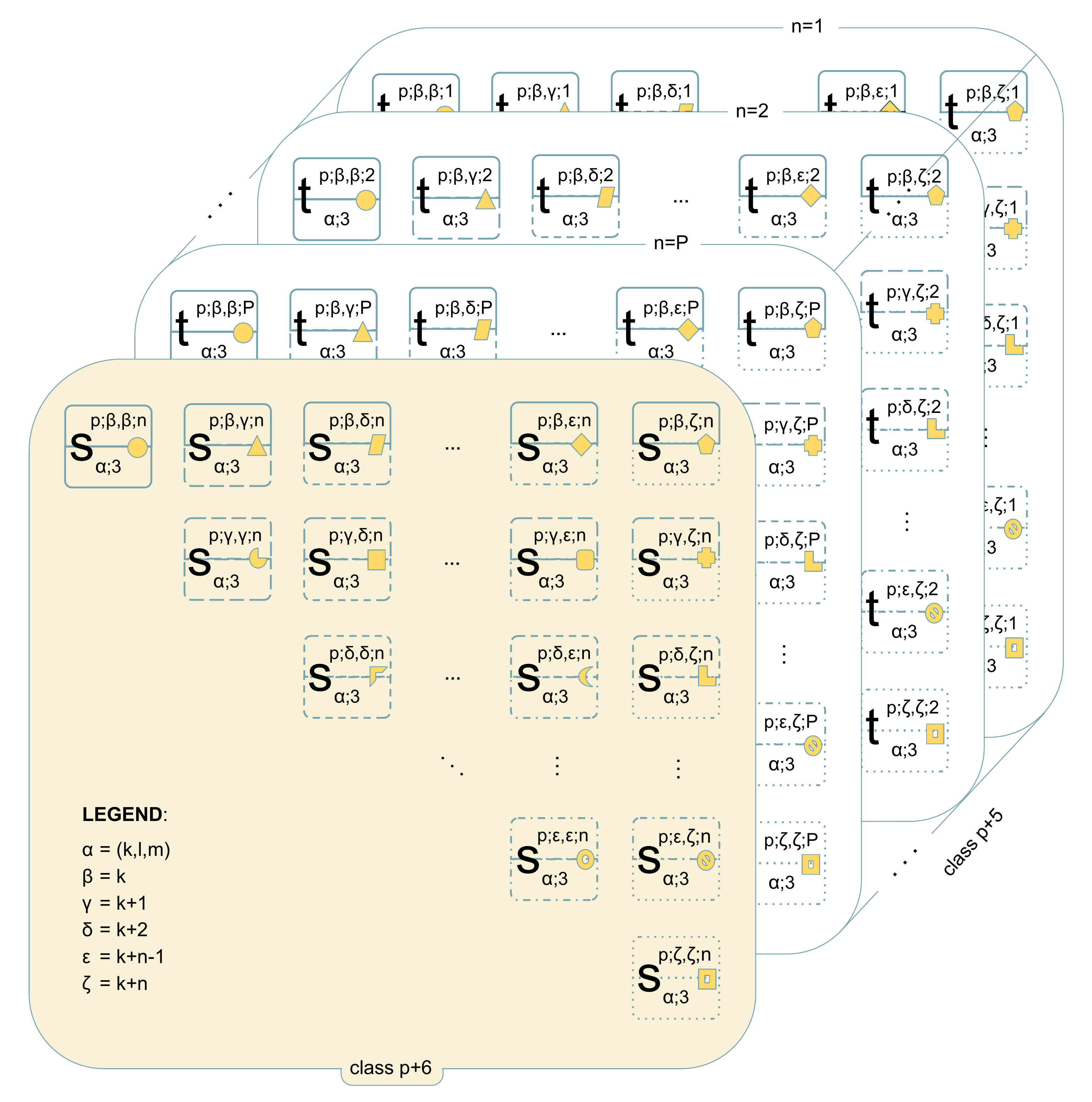}
    \caption{Relationships between classes $p+5$ and $p+6$ for $p$ - order functions.
    Each task from class $p+6$ corresponds to the approximation of the function value using Gaussian quadrature, therefore it depends on $P$ tasks from class $p+5$.
    The task $s$ depends on all tasks $t$ with regards to its distribution in subsequent sheets $1, 2, \dots, P$.}
    \label{fig:n_part7}
\end{figure}

Next, we construct two pairs of classes $p+3$ and $p+4$ (Figures \ref{fig:n_part4}, \ref{fig:n_part5}), and $p+5$ and $p+6$ (Figures \ref{fig:n_part6}, \ref{fig:n_part7}), in similar manner to classes $p+1$ and $p+2$.
Finally, we present all tasks in Tables \ref{table66} and \ref{table67}.

\begin{table}[ht!]
\centering
\begin{tabular}{|m{1.2cm}|m{7cm}|m{3.5cm}|}
    \hline
    $ s_{\alpha;3}^{p;\beta} $ & $ B_{\beta;p}(x^n) = B_{m;p}(x_{3}^{n}) \, B_{c;p}(x_{3}^{n})$ & \vtop{\hbox{}\hbox{}\hbox{$n \in \{1,2,\dots,P\}$,}\hbox{$\beta=\{k,l,m\} \in \mathcal{K}^\Delta_\alpha$} \hbox{$\gamma=\{a,b,c\} \in \mathcal{K}^\Delta_\alpha$}\hbox{$ I(\beta) \geq I(\gamma)$ }}  \\
    \hline
    $ t_{\alpha;3}^{p;\beta;n} $ & $ B_{\beta;p}(x^n) = B_{m;p}(x_{3}^{n}) \, B_{c;p}(x_{3}^{n})$ & \vtop{\hbox{}\hbox{}\hbox{$n \in \{1,2,\dots,P\}$,}\hbox{$\beta=\{k,l,m\} \in \mathcal{K}^\Delta_\alpha$} \hbox{$\gamma=\{a,b,c\} \in \mathcal{K}^\Delta_\alpha$}\hbox{$ I(\beta) \geq I(\gamma)$ }}  \\
    \hline
    $ s_{\alpha;2}^{p;\beta} $ & $ B_{\beta;p}(x^n) = B_{l;p}(x_{2}^{n}) \, B_{b;p}(x_{2}^{n})$ & \vtop{\hbox{}\hbox{}\hbox{$n \in \{1,2,\dots,P\}$,}\hbox{$\beta=\{k,l,m\} \in \mathcal{K}^\Delta_\alpha$} \hbox{$\gamma=\{a,b,c\} \in \mathcal{K}^\Delta_\alpha$}\hbox{$ I(\beta) \geq I(\gamma)$ }}  \\
    \hline
    $ t_{\alpha;2}^{p;\beta;n} $ & $ B_{\beta;p}(x^n) = B_{l;p}(x_{2}^{n}) \, B_{b;p}(x_{2}^{n})$ & \vtop{\hbox{}\hbox{}\hbox{$n \in \{1,2,\dots,P\}$,}\hbox{$\beta=\{k,l,m\} \in \mathcal{K}^\Delta_\alpha$} \hbox{$\gamma=\{a,b,c\} \in \mathcal{K}^\Delta_\alpha$}\hbox{$ I(\beta) \geq I(\gamma)$ }}  \\
    \hline
\end{tabular}
\caption{Computational tasks for performing computations of sum factorization algorithm of 3D order $p$ basis functions over element $E_{\alpha}$, $(k,l,m)=\alpha \in \mathcal{K}^\Delta$.
Part 1.}
\label{table66}
\end{table}    
\begin{table}[ht!]
\centering
\begin{tabular}{|m{1.2cm}|m{7cm}|m{3.5cm}|}
    \hline   
    $ s_{\alpha;1}^{p;\beta} $ & $ B_{\beta;p}(x^n) = B_{k;p}(x_{1}^{n}) \, B_{a;p}(x_{1}^{n})$ & \vtop{\hbox{}\hbox{}\hbox{$n \in \{1,2,\dots,P\}$,}\hbox{$\beta=\{k,l,m\} \in \mathcal{K}^\Delta_\alpha$} \hbox{$\gamma=\{a,b,c\} \in \mathcal{K}^\Delta_\alpha$}\hbox{$ I(\beta) \geq I(\gamma)$ }}  \\
    \hline
    $ t_{\alpha;1}^{p;\beta;n} $ & $ B_{\beta;p}(x^n) = B_{k;p}(x_{1}^{n}) \, B_{a;p}(x_{1}^{n})$ & \vtop{\hbox{}\hbox{}\hbox{$n \in \{1,2,\dots,P\}$,}\hbox{$\beta=\{k,l,m\} \in \mathcal{K}^\Delta_\alpha$} \hbox{$\gamma=\{a,b,c\} \in \mathcal{K}^\Delta_\alpha$}\hbox{$ I(\beta) \geq I(\gamma)$ }}  \\
    \hline 
    $ t_{\alpha;d}^{p;r;n} $ & $ B_{r;p}(x_{d}^{n}) $  & \vtop{\hbox{}\hbox{}\hbox{$d \in \{1,2,3\}$,}\hbox{$n \in \{1,2,\dots,P\}$ }}\\
    \hline
    \vdots&\vdots&\vdots\\
    \hline
    $ t_{\alpha;d}^{1;r;n} $ & $ B_{r;1}(x_{d}^{n}) $  & \vtop{\hbox{}\hbox{}\hbox{$d \in \{1,2,3\}$,}\hbox{$n \in \{1,2,\dots,P\}$ }}\\
    \hline
    $ t_{\alpha;d}^{0;r;n} $ & $ B_{r;0}(x_{d}^{n}) $ & \vtop{\hbox{}\hbox{}\hbox{$ \in \{1,2,3\}$,}\hbox{$n \in \{1,2,\dots,P\}$ }}\\
    \hline
\end{tabular}
\caption{Computational tasks for performing computations of sum factorization algorithm of 3D order $p$ basis functions over element $E_{\alpha}$, $(k,l,m)=\alpha \in \mathcal{K}^\Delta$.
Part2.}
\label{table67}
\end{table}

\subsection{Scheduling algorithm}
\label{sec:scheduling_algorithm}

To obtain a similar scheduling quality to the classical algorithm on massively parallel shared-memory machines, we employ the Foata-Normal-Form (FNF) \cite{traces}.
The Diekert dependency graphs  (see Section \ref{sec:porder}) show the consecutive Foata classes for each considered case of sum factorization.
Within a given Foata class, tasks can be executed in any order.
Completion of the entire previous Foata class is a sufficient condition to begin the computation of the next one.
The proposed strategy ensures no deadlocks, high-quality scheduling, and no need for intra-class synchronization.

Based on Figures \ref{fig:n_part1}, \ref{fig:n_part2} and \ref{fig:n_part3}, we can describe a general procedure for creating subsequent Foata classes, containing the following tasks:
\begin{itemize}
\item Class $m$, where $m \in \{0,\dots,p-1\}$
\begin{equation}
    \{  t_{\alpha;d}^{m;r;n} ;
    d\in \{1,2,3\}, \, n \in \{ 1,2,\dots,P\}  \},
    \label{classm}
\end{equation}
\item Class $p$
\begin{equation}
    \{  t_{\alpha;d}^{p;r;n} , s_{\alpha}^{p,n} ;
    d\in \{1,2,3\}, \, n \in \{ 1,2,\dots,P\}  \},
    \label{classp}
\end{equation}
\item Class $p+1$
\begin{equation}
    \{  t_{\alpha;1}^{p;\beta,\gamma;n};
    (\beta,\gamma) \in \mathcal{K}^\Delta_\alpha \times \mathcal{K}^\Delta_\alpha, \, I(\beta) \geq I(\gamma), \, n \in \{1,2,\dots,P\}  \},
    \label{classp1}
\end{equation}
\item Class $p+2$
\begin{equation}
    \{  s_{\alpha;1}^{p;\beta,\gamma};
    (\beta,\gamma) \in \mathcal{K}^\Delta_\alpha \times \mathcal{K}^\Delta_\alpha, \, I(\beta) \geq I(\gamma), \, n \in \{1,2,\dots,P\}  \},
    \label{classp2}
\end{equation}
\item Class $p+3$
\begin{equation}
    \{  t_{\alpha;2}^{p;\beta,\gamma;n};
    (\beta,\gamma) \in \mathcal{K}^\Delta_\alpha \times \mathcal{K}^\Delta_\alpha, \, I(\beta) \geq I(\gamma), \, n \in \{1,2,\dots,P\}  \},
    \label{classp3}
\end{equation}
\item Class $p+4$
\begin{equation}
    \{  s_{\alpha;2}^{p;\beta,\gamma};
    (\beta,\gamma) \in \mathcal{K}^\Delta_\alpha \times \mathcal{K}^\Delta_\alpha, \, I(\beta) \geq I(\gamma), \, n \in \{1,2,\dots,P\}  \},
    \label{classp4}
\end{equation}
\item Class $p+5$
\begin{equation}
    \{  t_{\alpha;3}^{p;\beta,\gamma;n};
    (\beta,\gamma) \in \mathcal{K}^\Delta_\alpha \times \mathcal{K}^\Delta_\alpha, \, I(\beta) \geq I(\gamma), \, n \in \{1,2,\dots,P\}  \},
    \label{classp5}
\end{equation}
\item Class $p+6$
\begin{equation}
    \{  s_{\alpha;3}^{p;\beta,\gamma};
    (\beta,\gamma) \in \mathcal{K}^\Delta_\alpha \times \mathcal{K}^\Delta_\alpha, \, I(\beta) \geq I(\gamma), \, n \in \{1,2,\dots,P\}  \},
    \label{classp6}
\end{equation}
\end{itemize}

The first Foata classes (\ref{classm}, \ref{classp1}) are responsible for valuating the values of 1D $n$-order basis functions over the element $E_\alpha$, at the Gaussian quadrature points, using recursive Cox--de--Boor formulae (\ref{eq:cox1}, \ref{eq:cox2}) and the Jacobian (\ref{classp1}).
Subsequent Foata classes of two kinds follow this:
\begin{enumerate}
    \item Computational tasks (\ref{classp1}, \ref{classp3}, \ref{classp5}) evaluating the values of the dot products of 1D $n$-order basis functions over the element $E_\alpha$ at Gaussian quadrature point,
    \item computational tasks (\ref{classp2}, \ref{classp4}, \ref{classp6}) evaluating the values of the dot products of 1D $n$-order basis functions over the element $E_\alpha$ and buffers.
\end{enumerate}

All the tasks mentioned above are performed on a homogeneous architecture.
Thus, we can expect near-identical execution time for each of them inside a particular Foata class.
Consequently, all tasks from the particular Foata class can be effectively scheduled as a common bag.

Over each element $E_\alpha$ we repeat the same procedure of invoking tasks using parameters associated with this element.
We invoke Foata classes starting from the Foata class 0, and each time wait for all tasks to be completed before invoking the next Foata class.
Using a simplified scheduling method, based on FNF and the proposed above, despite having no theoretical proof, results in near-optimal performance in practical applications while maintaining a relatively simple implementation.

\section{Numerical results}
\label{sec:numres}

Now, we compare the computational performance of parallel integration using the classical algorithm and sum factorization.
In both cases, implementation was done in Fortran 2003, using OpenMP for loop parallelization.
The measurements concern the execution time for the sequential integration algorithm executed on CPU and the concurrent integration algorithm run on a shared memory CPU with 12 cores.
Computations were performed on a Banach Linux workstation equipped with AMD Ryzen 9 3900X processor and 64GB RAM.
It is worth noting that the CPU, despite having 3.8 GHz base clock speed and 4.6 GHz boost, was working at a constant 4.0 GHz in the multi-threaded (12 cores) workload and at 4.1 GHz in single-threaded workload (1 core).
The computations have been performed using the code compiled with ifort with -O2 level of optimization.

In Sections \ref{sec:inside_element}, \ref{sec:over_element}, and \ref{sec:amdahl} we present the experimental results.
In Section \ref{sec:discussion} we discuss obtained results.

\subsection{Inside element scalability}
\label{sec:inside_element}

\begin{figure}[ht!]
\centering
\includegraphics[width=0.95\textwidth]{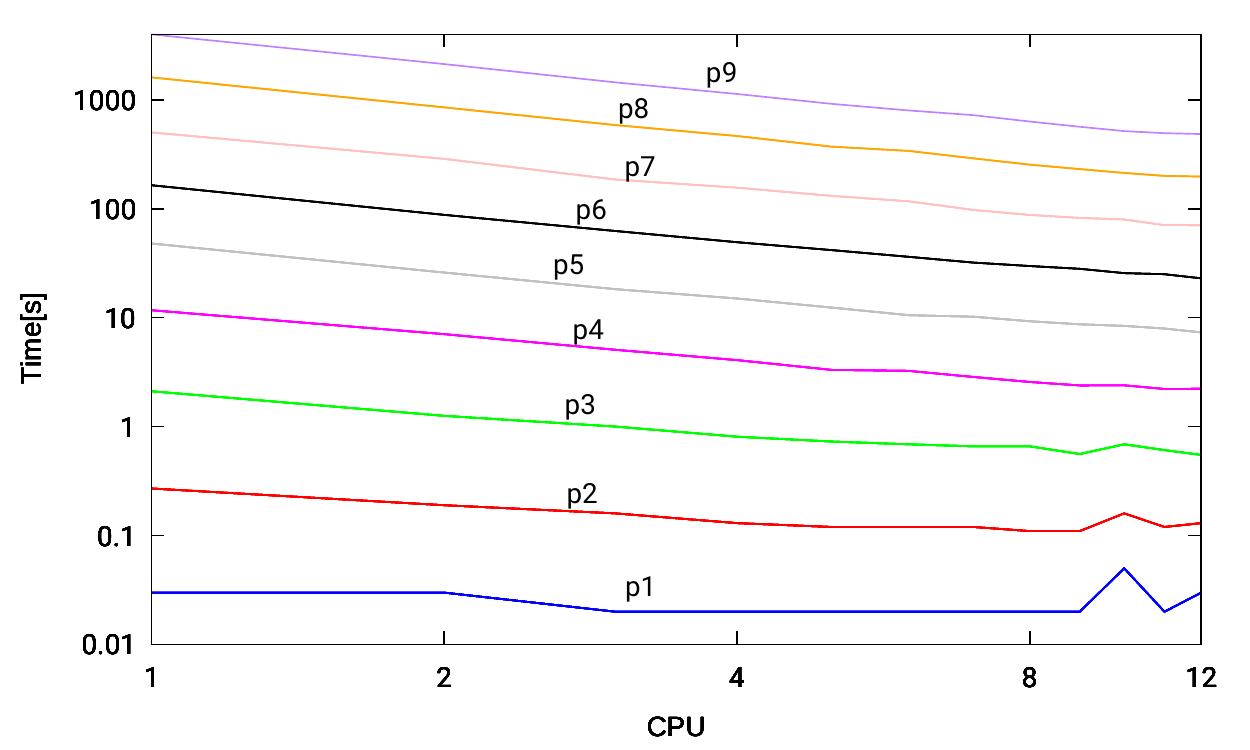}
\caption{Strong scaling time for classical integration algorithm.
Computations performed on $20^3$ elements mesh, different polynomial orders.
Parallelism inside element.}
\label{fig:scal_class}
\end{figure}

\begin{figure}[ht!]
\centering
\includegraphics[width=0.95\textwidth]{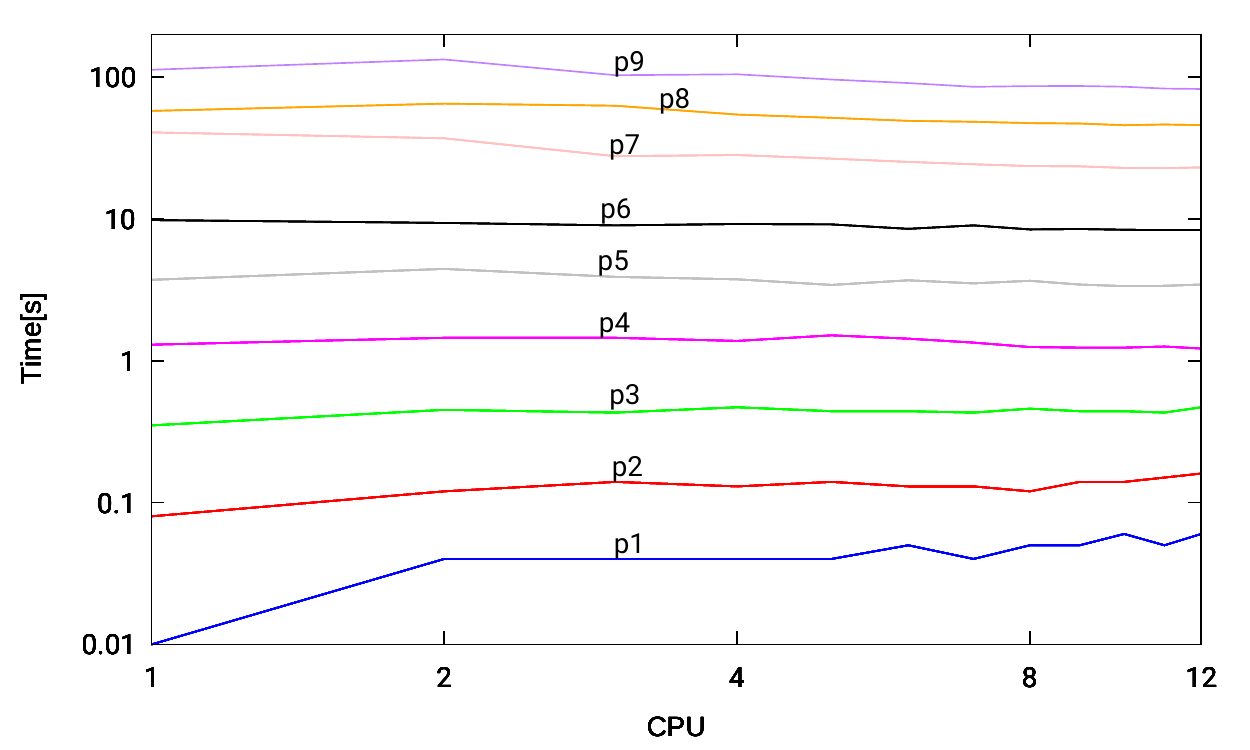}
\caption{Strong scaling time for sum factorization algorithm.
Computations performed on $20^3$ elements mesh, different polynomial orders.
Parallelism inside element.}
\label{fig:scal_sum}
\end{figure}

\begin{figure}[ht!]
\centering
\includegraphics[width=0.95\textwidth]{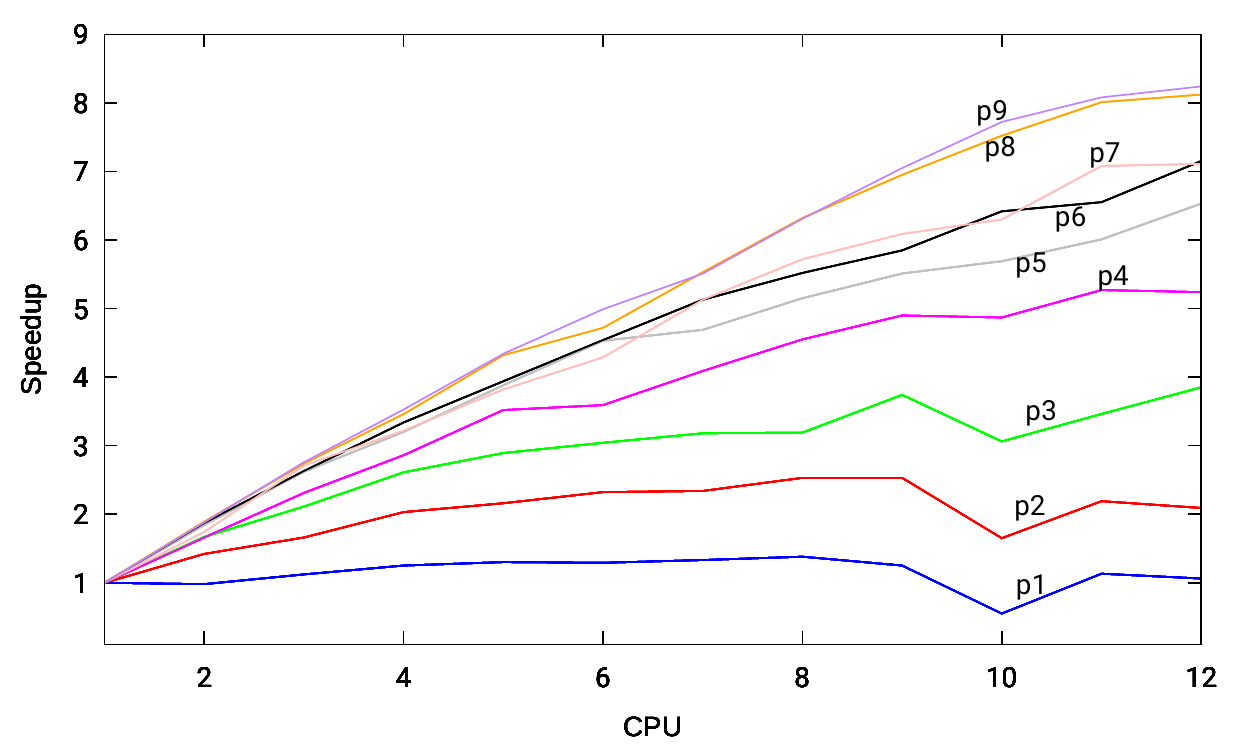}
\caption{Strong scaling speedup for classical integration algorithm.
Computations performed on $20^3$ elements mesh, different polynomial orders.
Parallelism inside element.}
\label{fig:spd_class}
\end{figure}

\begin{figure}[ht!]
\centering
\includegraphics[width=0.95\textwidth]{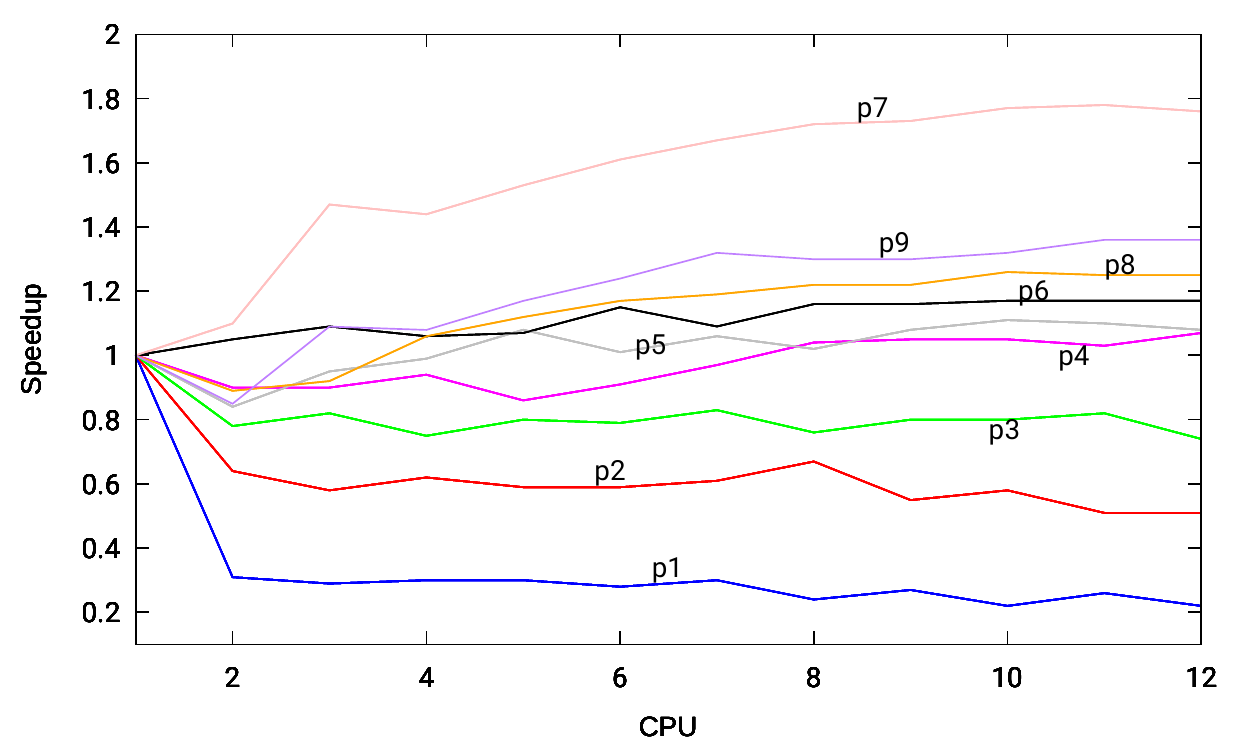}
\caption{Strong scaling speedup for sum factorization algorithm.
Computations performed on $20^3$ elements mesh, different polynomial orders.
Parallelism inside element.}
\label{fig:spd_sum}
\end{figure}

\begin{figure}[ht!]
\centering
\includegraphics[width=0.95\textwidth]{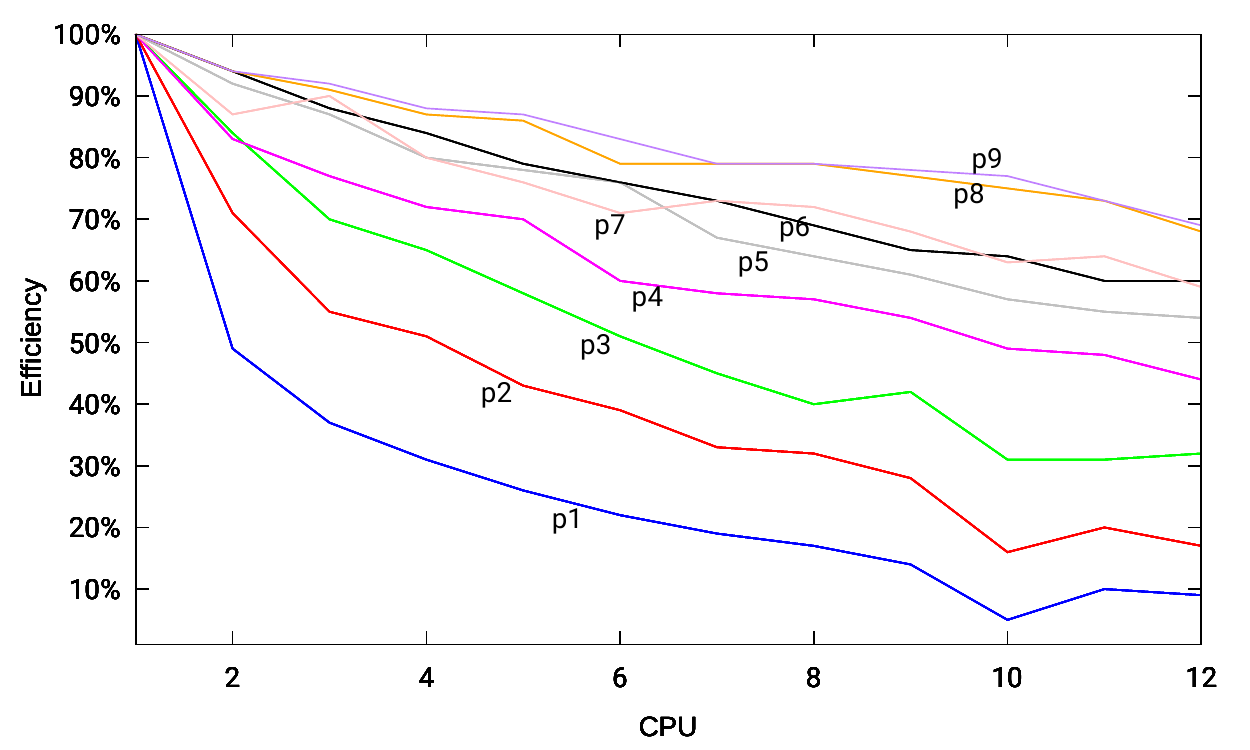}
\caption{Strong scaling efficiency for classical integration algorithm.
Computations performed on $20^3$ elements mesh, different polynomial orders.
Parallelism inside element.}
\label{fig:eff_class}
\end{figure}

\begin{figure}[ht!]
\centering
\includegraphics[width=0.95\textwidth]{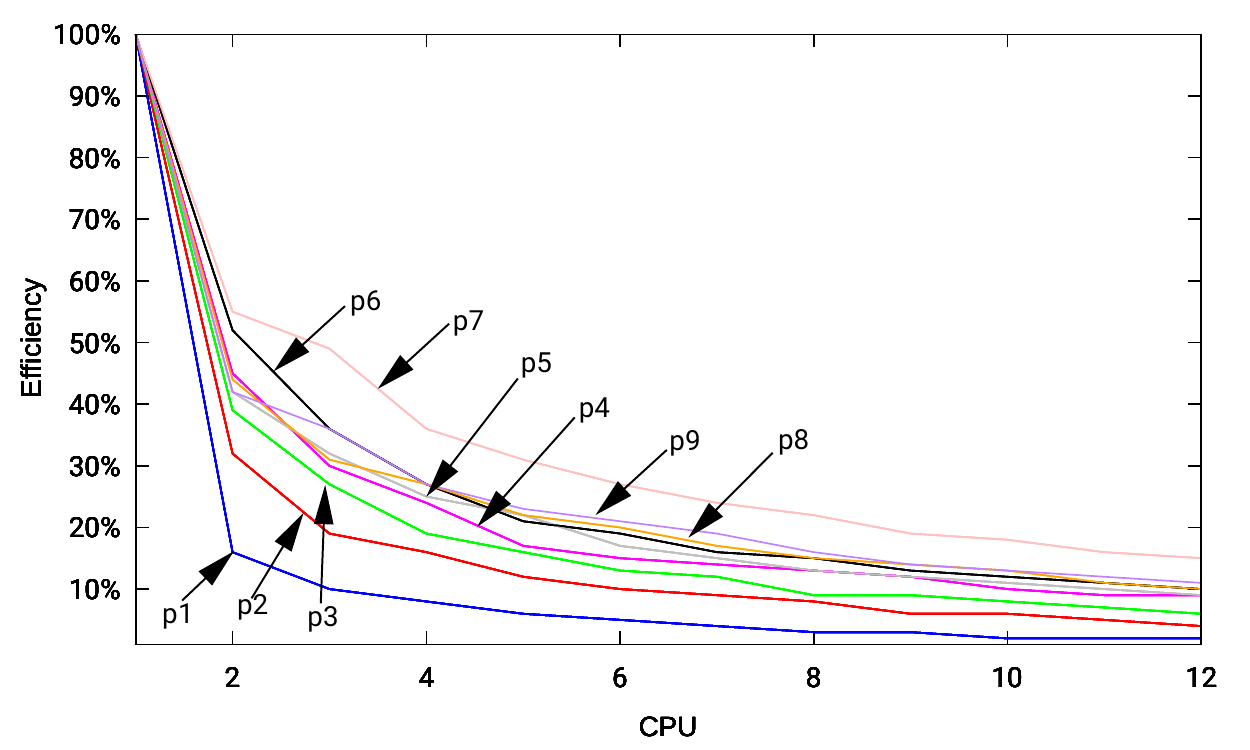}
\caption{Strong scaling efficiency for sum factorization algorithm.
Computations performed on $20^3$ elements mesh, different polynomial orders.
Parallelism inside element.}
\label{fig:eff_sum}
\end{figure}

We first performed computations with parallelization inside an element, then sequential looping over elements.
In such a case, we consider a mesh of $20^3$ elements.
The comparison of the scalability for different polynomial orders is presented in Figures \ref{fig:scal_class} and \ref{fig:scal_sum}.
Figures \ref{fig:spd_class} and \ref{fig:spd_sum} represent speedup.
Finally, in Figures \ref{fig:eff_class} and \ref{fig:eff_sum} we presented efficiency for the classical integration algorithm and sum factorization respectively.

\newpage

\subsection{Over element scalability}
\label{sec:over_element}

As a second experiment, we performed computations with sequential computations inside the element and parallel looping over elements.
In this case, we also used a mesh of $30^3$ elements.

\begin{figure}[ht!]
\centering
\includegraphics[width=0.95\textwidth]{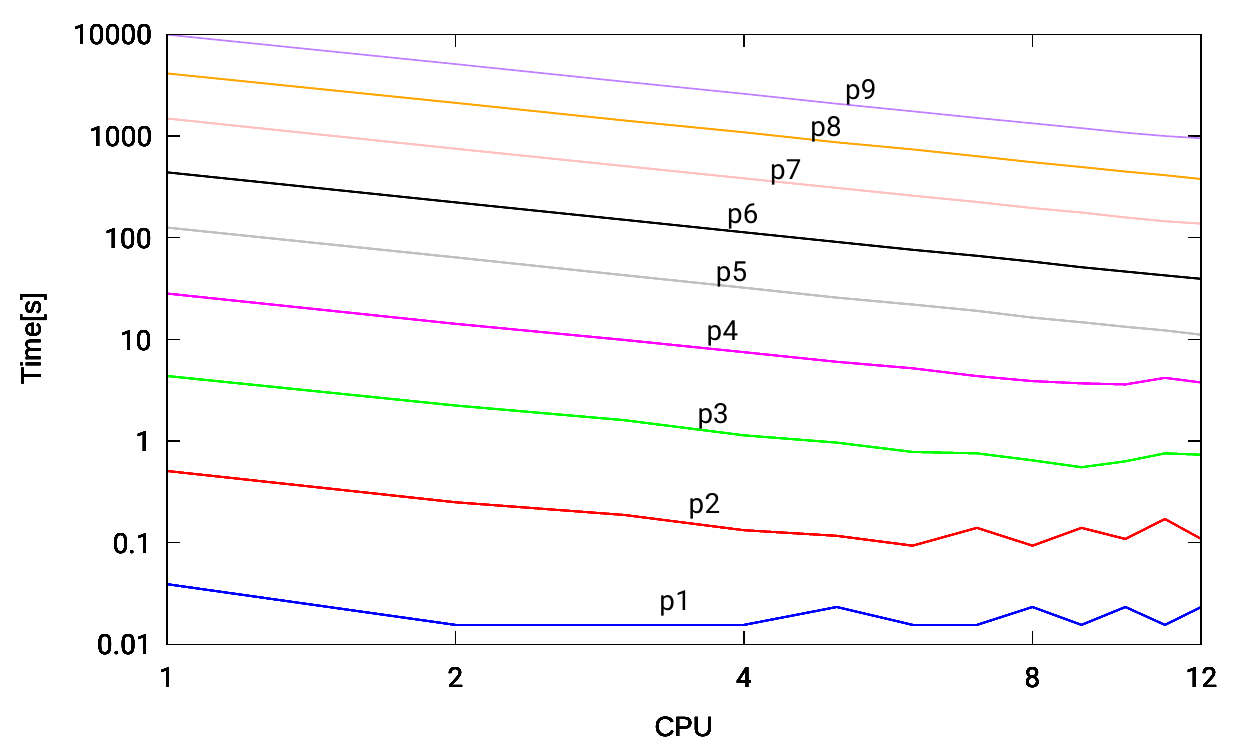}
\caption{Strong scaling time for classical integration algorithm.
Computations performed on $30^3$ elements mesh, different polynomial orders.
Parallelism over all elements.}
\label{fig:scal_class_elem}
\end{figure}

\begin{figure}[ht!]
\centering
\includegraphics[width=0.95\textwidth]{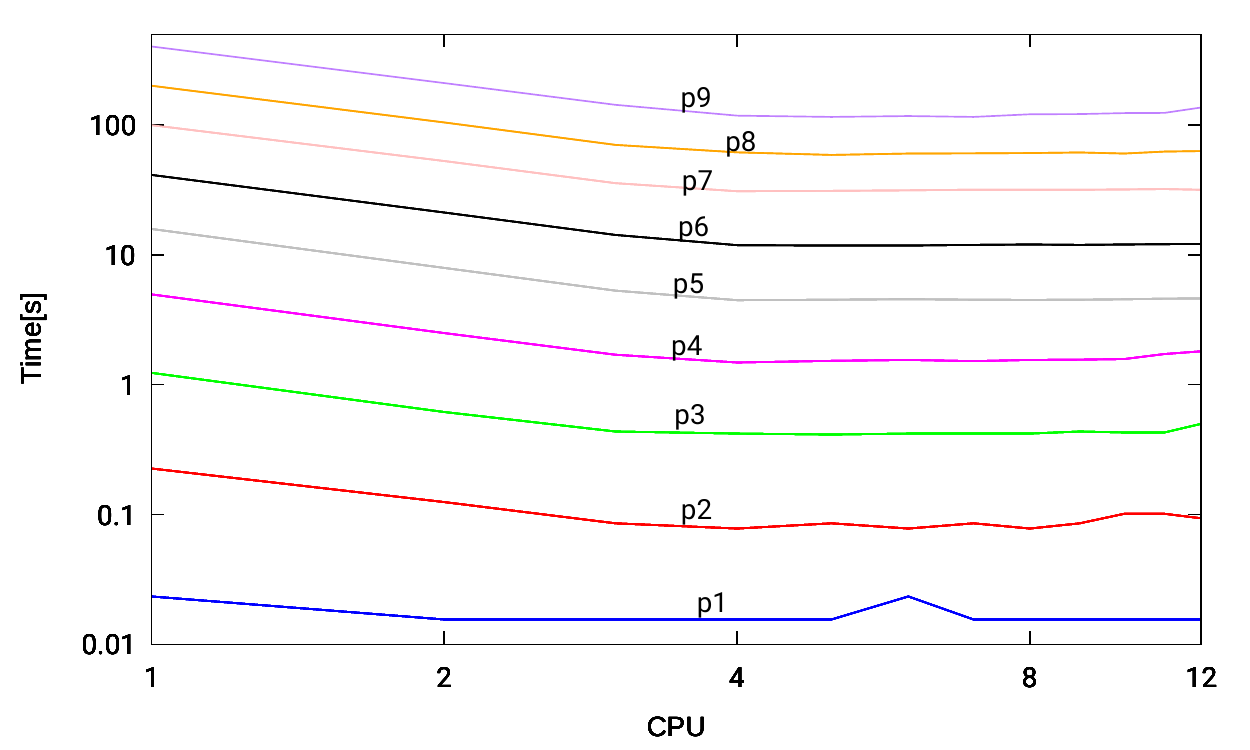}
\caption{Strong scaling time for sum factorization algorithm.
Computations performed on $30^3$ elements mesh, different polynomial orders.
Parallelism over all elements.}
\label{fig:scal_sum_elem}
\end{figure}

\begin{figure}[ht!]
\centering
\includegraphics[width=0.95\textwidth]{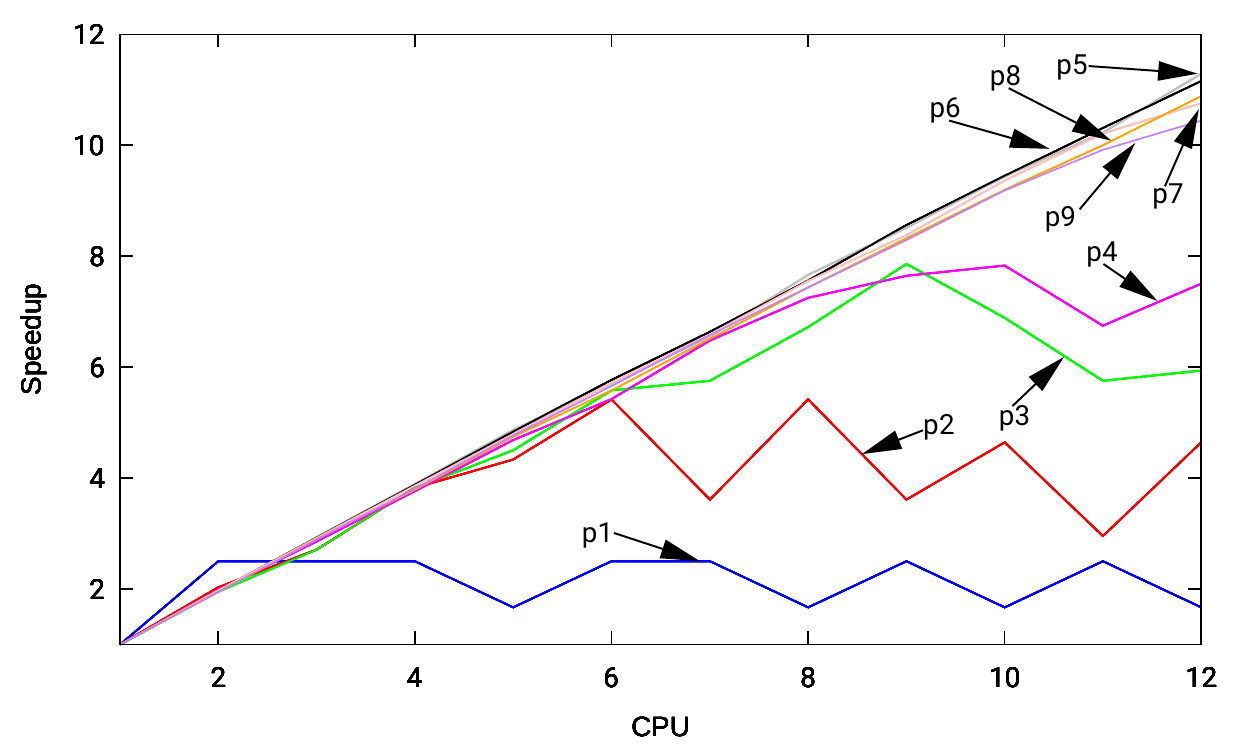}
\caption{Strong scaling speedup for classical integration algorithm.
Computations performed on $30^3$ elements mesh, different polynomial orders.
Parallelism over all elements.}
\label{fig:spd_class_elem}
\end{figure}

\begin{figure}[ht!]
\centering
\includegraphics[width=0.95\textwidth]{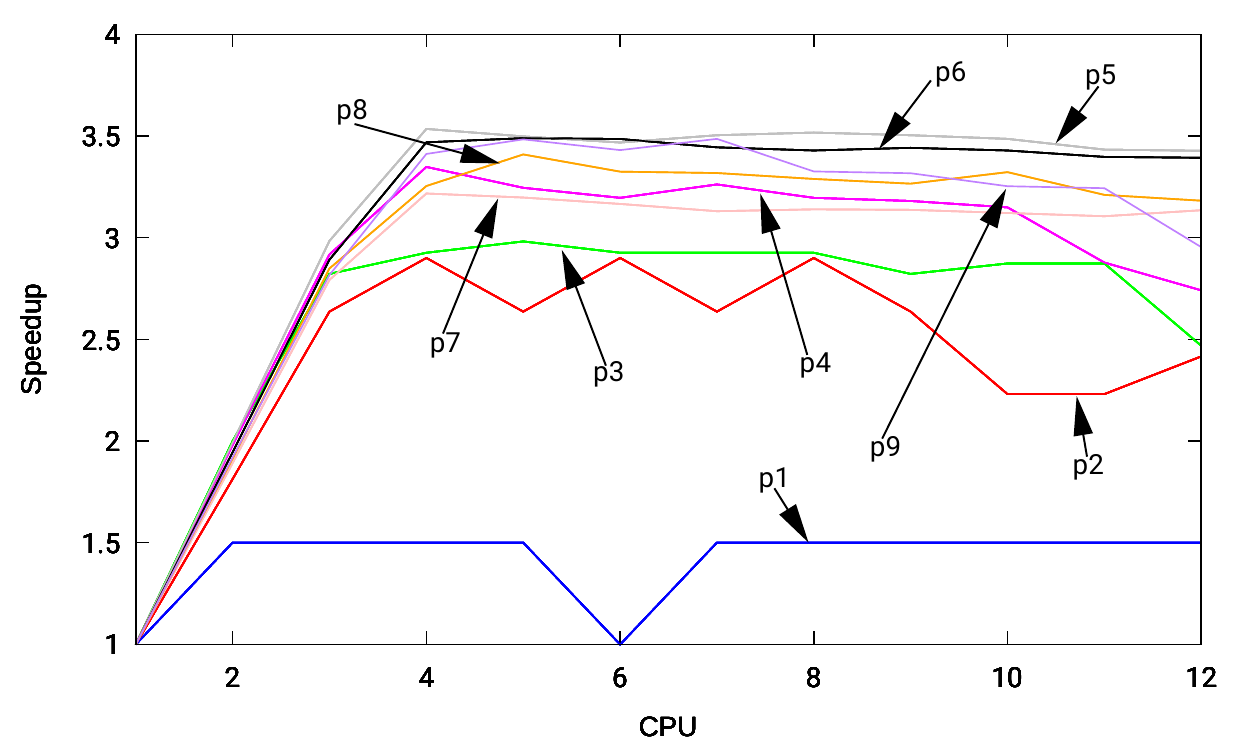}
\caption{Strong scaling speedup for sum factorization algorithm.
Computations performed on $30^3$ elements mesh, different polynomial orders.
Parallelism over all elements.}
\label{fig:spd_sum_elem}
\end{figure}

\begin{figure}[ht!]
\centering
\includegraphics[width=0.95\textwidth]{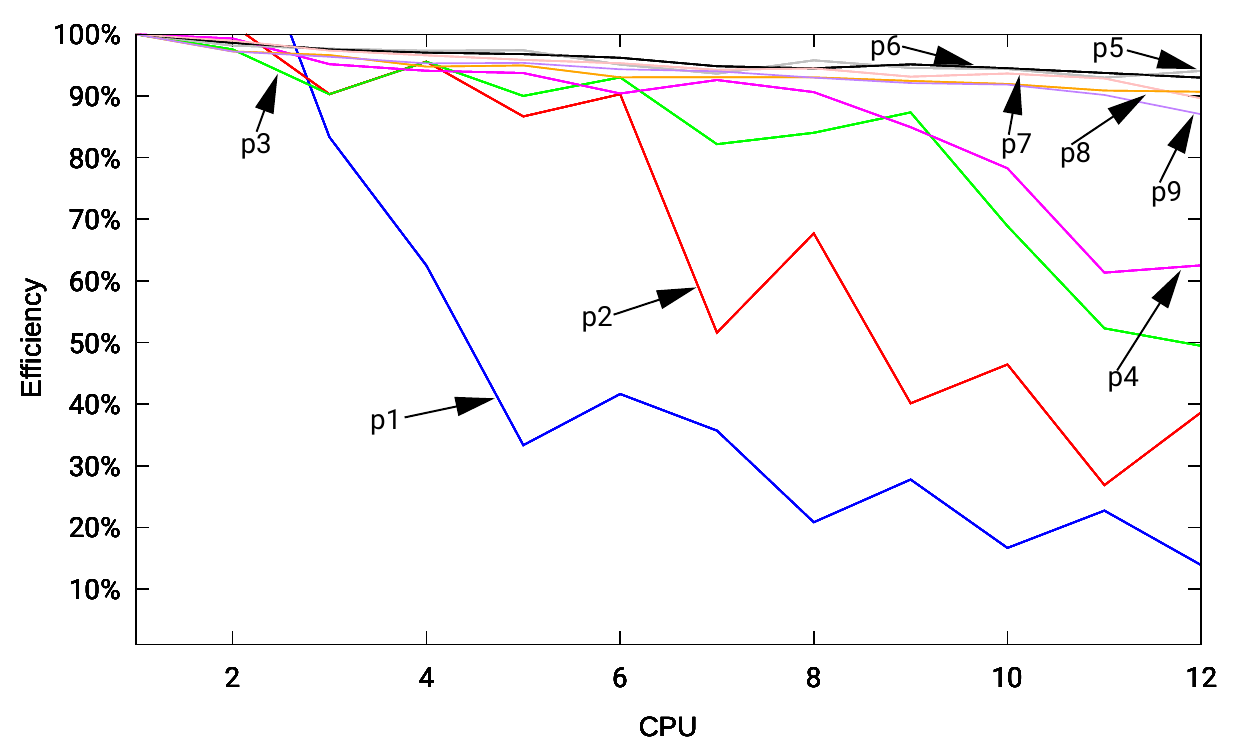}
\caption{Strong scaling efficiency for classical integration algorithm.
Computations performed on $30^3$ elements mesh, different polynomial orders.
Parallelism over all elements.}
\label{fig:eff_class_elem}
\end{figure}

\begin{figure}[ht!]
\centering
\includegraphics[width=0.95\textwidth]{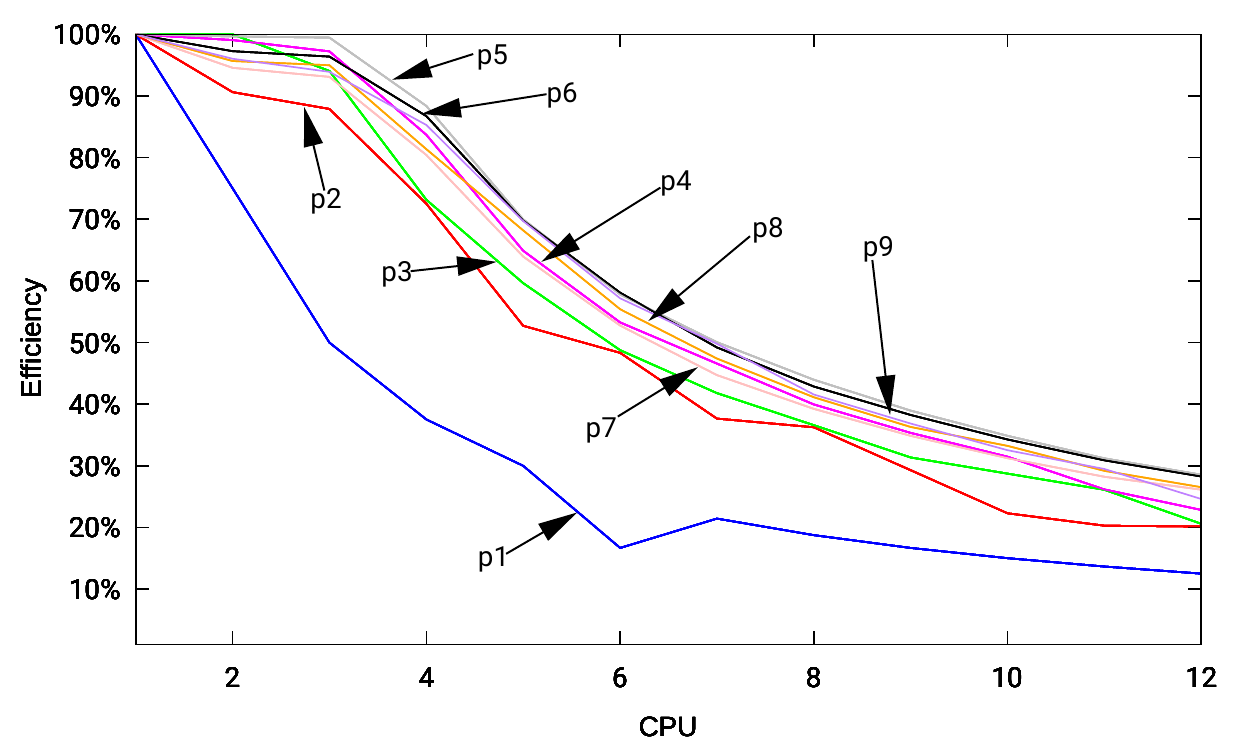}
\caption{Strong scaling efficiency for sum factorization algorithm.
Computations performed on $30^3$ elements mesh, different polynomial orders.
Parallelism over all elements.}
\label{fig:eff_sum_elem}
\end{figure}

The comparison of scaling for different polynomial orders is presented in Figures \ref{fig:scal_class_elem} and \ref{fig:scal_sum_elem}.
Figures \ref{fig:spd_class_elem} and \ref{fig:spd_sum_elem} represent the speedup.
Finally, in Figures \ref{fig:eff_class_elem} and \ref{fig:eff_sum_elem} we present the efficiency for the classical integration algorithm and sum factorization respectively.

\newpage

\subsection{Speedup limits}
\label{sec:amdahl}

As a final experiment, we estimate the maximum speedup for both the parallelization schemes (see Sections \ref{sec:inside_element} and \ref{sec:over_element}), and also its combination.
When considering integration inside a single element, the problem size is fixed regardless of the mesh size.
Amdahl's law is appropriate for this kind of scenario.
Therefore, to find the percentage of the algorithm which benefits from speedup $\mathcal{P}$, we invoke the Amdahl's equation:
\begin{equation}
    \mathcal{S}(\nu) = \frac{1}{ (1-\mathcal{P}) + \frac{\mathcal{P}}{\nu} }
    \label{eq:amdahl}
\end{equation}
where $\mathcal{P}$ denotes the percentage of the algorithm which benefits from the parallel speedup,
$\nu$ is the number of threads, and $\mathcal{S}(\nu)$ is the measured speedup when using $\nu$ threads.

From the previous equation, we can derive the value of $P$ and the speedup limit, which are explicitly given by:
\begin{equation}
    \mathcal{P} = \frac{ \frac{\nu}{\mathcal{S}(\nu)} -\nu }{ 1 - \nu},
    \label{eq:p-compute}
\end{equation}
and
\begin{equation}
    \mathcal{S}(\infty) = \lim\limits_{\nu \rightarrow \infty} \frac{1}{ (1-\mathcal{P}) + \frac{\mathcal{P}}{\nu}} = \frac{1}{1-\mathcal{P}},
    \label{eq:limit}
\end{equation}
respectively.
For different values of $p$, we consider the maximum experimental speedup observed from the numerical results for both methods.
Next, using equations (\ref{eq:p-compute}) and (\ref{eq:limit}), we computed the percentage of algorithm that benefits from the parallel speedup and the theoretical maximum speedup.
Finally, we estimated the combined maximum speedup by assuming two layers of parallelism.
This is, one layer representing the scheme of Section~\ref{sec:inside_element}, and another representing the scheme of Section~\ref{sec:over_element}.
The results for the classical integration algorithm are presented in Table~\ref{tab:speedup_classic}, while the results for sum factorization in Table~\ref{tab:speedup_sumfact}.

\begin{table}[ht!]
    \centering
    \begin{tabular}{c||c|c|c|c||c|c|c|c||c}
        $p$ &
        $\nu_i$ & $\mathcal{S}_{i}(\nu)$ & $\mathcal{P}_{i}$ & $\mathcal{S}_{i}(\infty)$ &
        $\nu_e$ & $\mathcal{S}_{e}(\nu)$ & $\mathcal{P}_{e}$ & $\mathcal{S}_{e}(\infty)$ & 
        $\mathcal{S}_{c}(\infty)$\\
        \hline
        1	& 8	    & 1.38	& 0.31	& 1.46	& 3	    & 2.5	& 0.9	& 10.00	    & 14.59 \\
        2	& 8	    & 2.53	& 0.69	& 3.24	& 6	    & 5.4	& 0.98	& 45.00 	& 145.69 \\
        3	& 12	& 3.85	& 0.81	& 5.20	& 9	    & 7.8	& 0.98	& 52.00	    & 270.21 \\
        4	& 11	& 5.27	& 0.89	& 9.20	& 10	& 7.8	& 0.97	& 31.91	    & 293.47 \\
        5	& 12	& 6.53	& 0.92	& 13.13	& 12	& 11.29	& 0.99	& 174.92	& 2296.93 \\
        6	& 12	& 7.15	& 0.94	& 16.22	& 12	& 11.15	& 0.99	& 144.29	& 2339.94 \\
        7	& 12	& 7.11	& 0.94	& 15.99	& 12	& 10.75	& 0.99	& 94.60	    & 1513.02 \\
        8	& 12	& 8.12	& 0.96	& 23.02	& 12	& 10.88	& 0.99	& 106.86	& 2459.92 \\
        9	& 12	& 8.24	& 0.96	& 24.11	& 12	& 10.44	& 0.99	& 73.62	    & 1774.60 \\
    \end{tabular}
    \caption{classical integration method.
    Bottom index $i$ stands for "inside element", $e$ over all elements, and $c$ combined.}
    \label{tab:speedup_classic}
\end{table}

\begin{table}[ht!]
    \centering
    \begin{tabular}{c||c|c|c|c||c|c|c|c||c}
        $p$ &
        $\nu_i$ & $\mathcal{S}_{i}(\nu)$ & $\mathcal{P}_{i}$ & $\mathcal{S}_{i}(\infty)$ &
        $\nu_e$ & $\mathcal{S}_{e}(\nu)$ & $\mathcal{P}_{e}$ & $S_{e}(\infty)$ & 
        $\mathcal{S}_{c}(\infty)$\\
        \hline
        1	& 1	    & 1	    & 0 	& 1 	& 2	& 1.5	& 0.67	& 3	    & 3.00 \\
        2	& 1	    & 1	    & 0	    & 1 	& 4	& 2.9	& 0.87	& 7.91	& 7.91 \\
        3	& 1	    & 1	    & 0	    & 1	    & 4	& 2.9	& 0.87	& 7.91	& 7.91 \\
        4	& 12	& 1.07	& 0.07	& 1.08	& 4	& 3.3	& 0.93	& 14.14	& 15.23 \\
        5	& 10	& 1.11	& 0.11	& 1.12	& 4	& 3.5	& 0.95	& 21	& 23.60 \\
        6	& 10	& 1.17	& 0.16	& 1.19	& 4	& 3.5	& 0.95	& 21	& 25.04 \\
        7	& 12	& 1.77	& 0.47	& 1.9	& 4	& 3.2	& 0.92	& 12	& 22.84 \\
        8	& 10	& 1.26	& 0.23	& 1.3	& 4	& 3.4	& 0.94	& 17	& 22.06 \\
        9	& 11	& 1.36	& 0.29	& 1.41	& 4	& 3.5	& 0.95	& 21	& 29.63 \\
    \end{tabular}
    \caption{Sum factorization.
    Bottom index $i$ stands for "inside element", $e$ over all elements, and $c$ combined.}
    \label{tab:speedup_sumfact}
\end{table}

\newpage

\subsection{Discussion of the numerical results}
\label{sec:discussion}

For different values of $p$, we consider the maximum experimental speedup observed from the numerical results for both methods.
Next, using equations (\ref{eq:p-compute}) and (\ref{eq:limit}), we computed the percentage of algorithm that benefits from the parallel speedup and the theoretical maximum speedup.
From Figures \ref{fig:spd_class} and \ref{fig:spd_class_elem} we can observe outstanding speedup for classical method in both scenarios of parallelism.
Furthermore Figures \ref{fig:eff_class} and \ref{fig:eff_class_elem} proven high efficiency of hardware utilization.
Figures present increased parallel performance (speedup and efficiency) for higher polynomial order ($p$) B-spline basis functions.

Figures \ref{fig:scal_sum} and \ref{fig:spd_sum} present unexpected behaviour of sum factorization with parallel loops inside elements.
Even parallel loops over all elements, presented in Figure \ref{fig:spd_sum_elem} scale up to 4 cores with expected behavior. Above four cores, speedup remains at a constant level.
This corresponds with low efficiency in multicore applications, as can be seen in Figures \ref{fig:eff_sum} and \ref{fig:eff_sum_elem}.

From Tables \ref{tab:speedup_classic} and \ref{tab:speedup_sumfact}, we can observe that the theoretical maximum speedup for the classical method behaves similarly to the results presented in  \cite{parallel_integration}.
In Diekert graphs (Figures \ref{fig:n_part1}-\ref{fig:n_part7}), it can be observed that sum factorization requires a multitude more memory synchronizations than the classical method.

We also compare computational times for the classical integration and the sum factorization in several scenarios.
We focused on $p=9$ since, theoretically, it should be the best scenario of sum factorization.
We take into consideration three scenarios for a $30^3$ mesh size;
1) Single-core CPU execution,
2) Shared memory CPU computations,
3) (Multiple) GPU execution.
Classical integration on single-core takes 9931.758 seconds,
12 core OpenMP implementation takes 951 seconds,
and estimated GPU implementation should take 4.596 seconds.
Sum factorization integration on a single core takes 403.586 seconds,
Four-core OpenMP implementation takes 118.296 seconds
and estimated GPU implementation should take 13.62 seconds.

\section{Conclusions}
\label{sec:conclussions} 
 
In terms of computational performance, we discussed and compared two standard methods used for the integration in IGA-FEM; the classical integration method and sum factorization.
For the comparison, we considered several scenarios of performing a shared memory layer of computations on hybrid memory clusters.
First, we consider a single-core implementation as the baseline.
Then, we measure experimental performance in two ways of parallel integration in shared memory, using OpenMP, with parallel loops over elements and parallel loops inside elements.
In the final scenario, we estimate performance on massively parallel shared-memory machines, such as GPU, by combining maximum scalability estimates (see Section \ref{sec:amdahl}).

As expected, when assigned to a specific computational node, the sum factorization method performs better than the classical integration method.
From the numerical results with a polynomial degree $p=9$, being the worst-case scenario from the considered experiments, we can observe that the classical method is approximately 70 times slower than the sum factorization method in both scenarios of parallel integration in shared memory.
Even though, when comparing single-core sum factorization with parallelized on 12 CPU cores classical integration method, still sum factorization is the clear winner.

When considering parallelized loops inside the elements, we observe very efficient parallelization for the classical integration method.
However, sum factorization does not parallelize as expected.
Indeed, we observe an evident loss in performance when considering more than one core.
Additionally, when considering the standard loops over elements, we observe performance gain for sum factorization only up to 4 cores in a shared memory (see Figures \ref{fig:spd_sum}, \ref{fig:spd_sum_elem}).

Finally, based on the previous work \cite{parallel_integration}, we can assume that estimate the performance for both parallelization methods mixed on massively parallel machines, such as GPUs.
In such a case, the classical integration method parallelizes outstandingly, resulting in faster execution than sum factorization.
In other words, numerical results show that the classical integration method running on a GPU can be faster than sum factorization by one or two orders of magnitude.
A possible explanation for this small performance gain, or lack of such in some cases for sum factorization, is possibly limited by the memory synchronization and the memory access.
Despite the higher computational cost of the classical method concerning sum factorization, such a method requires fewer data dependencies and synchronizations than sum factorization.
However, when considering low cores machines, sum factorization is the method of choice over the classical one.
The best parallelization strategy we observe in such a case is to use 4 CPU cores in shared memory.

\vspace{0.5cm}
\noindent{\bf{Acknowledgments}}
This project has received funding from the European Union's Horizon 2020 research and innovation programme under the Marie Sklodowska-Curie grant agreement No 777778 (MATHROCKS).
The work of SR has also been partially supported by the Chilean grant ANID Fondecyt No 3210009.


\bibliographystyle{elsarticle-num}
\bibliography{bibliography}

\end{document}